\documentclass{article}
\usepackage{hyperref}
\usepackage{amsmath,amsthm}
\usepackage[dvips]{graphicx} % Para a inclusão de gráficos no formato EPS
\graphicspath{{eps/}}
\usepackage{amssymb}
\usepackage[latin1]{inputenc}

\newtheorem{theorem}{Theorem}

\newtheorem*{definition}{Definition}
\newtheorem{remark}[theorem]{Remark}

\title{Splitting of separatrices for the Hamiltonian-Hopf bifurcation
with the Swift-Hohenberg equation as an example}

\author{José Pedro Gaivão\thanks{Supported by FCT - Fundação para a Ciência e Tecnologia, Portugal, with grant SFRH/BD/30596/2006.}\: and Vassili Gelfreich\\
\small Mathematics Institute, University of Warwick,
\small Coventry, CV4 7AL, UK\\
\small \begin{tabular}{ll}
E-mail:&J.P.Romana-Gaivao@warwick.ac.uk\\ &V.Gelfreich@warwick.ac.uk
\end{tabular}}
\begin{document}
\maketitle
\begin{abstract}
We study homoclinic orbits of the Swift-Hohenberg equation near a Hamiltonian-Hopf bifurcation.
It is well known that in this case the normal form of the equation is integrable at all orders.
Therefore the difference between the stable and unstable manifolds is exponentially small and 
the study requires a method capable to detect phenomena beyond all algebraic orders provided by the normal form theory.
We propose an asymptotic expansion for an homoclinic invariant which quantitatively describes the
transversality of the invariant manifolds. We perform high-precision numerical experiments to support  
validity of the asymptotic expansion and evaluate a Stokes constant numerically using two independent methods. 
\end{abstract}

\section{The generalized Swift-Hohenberg equation}

The generalized Swift-Hohenberg equation (GSHE),
\begin{equation}\label{4:GSHE}
u_t = \epsilon u + \kappa u^2-u^3-(1+\Delta)^2 u 
\end{equation}
is widely used to model nonlinear phenomena 
in various areas of modern Physics including hydrodynamics, pattern formation 
and nonlinear optics (e.g. \cite{BurkeK2006,HaragusS2007}). 
This equation (with $\kappa=0$) was originally introduced 
 by Swift and Hohenberg~\cite{JSPH:77}
in a study of thermal fluctuations in a convective instability. 

\begin{figure}[t]
  \begin{center}
    \includegraphics[width=4in]{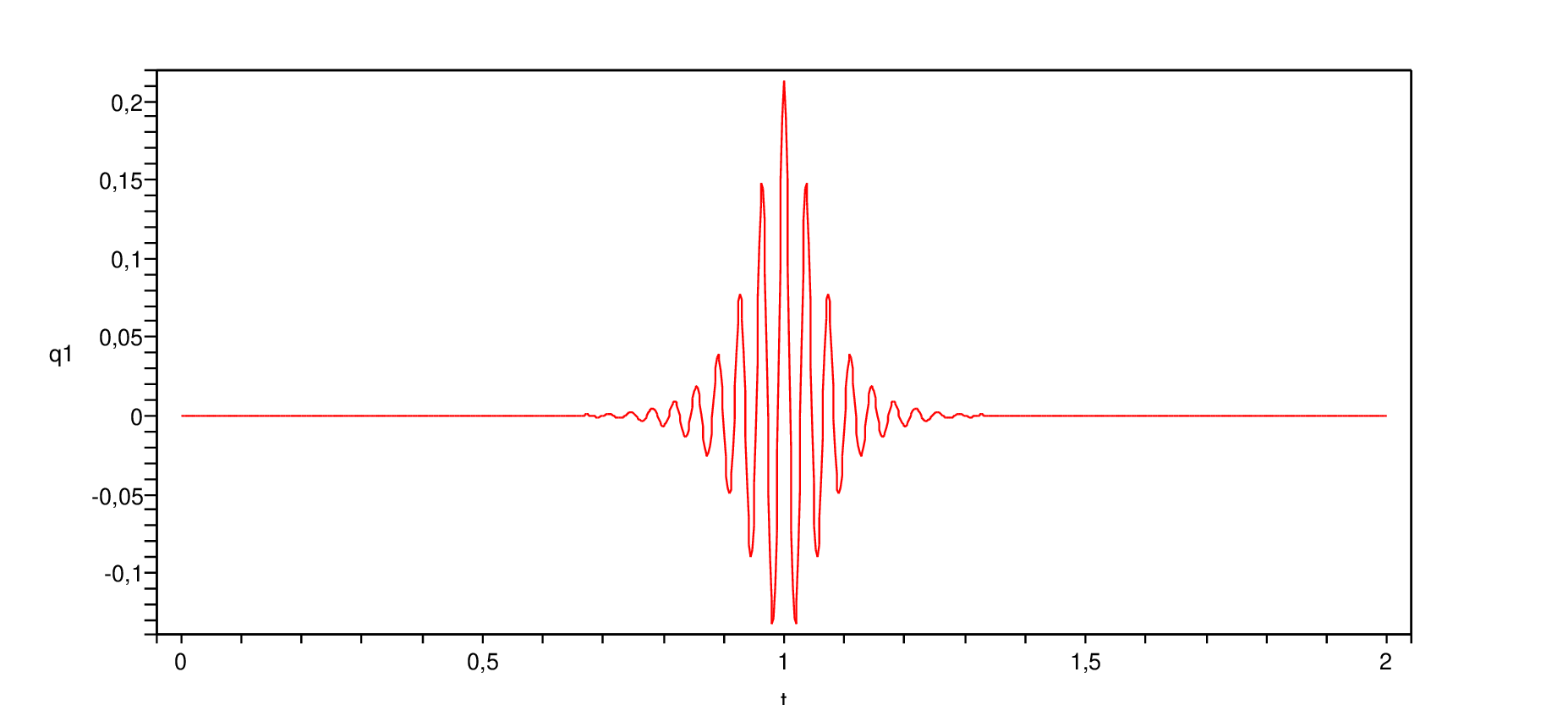}
  \end{center}
\begin{center}
    \includegraphics[width=4in]{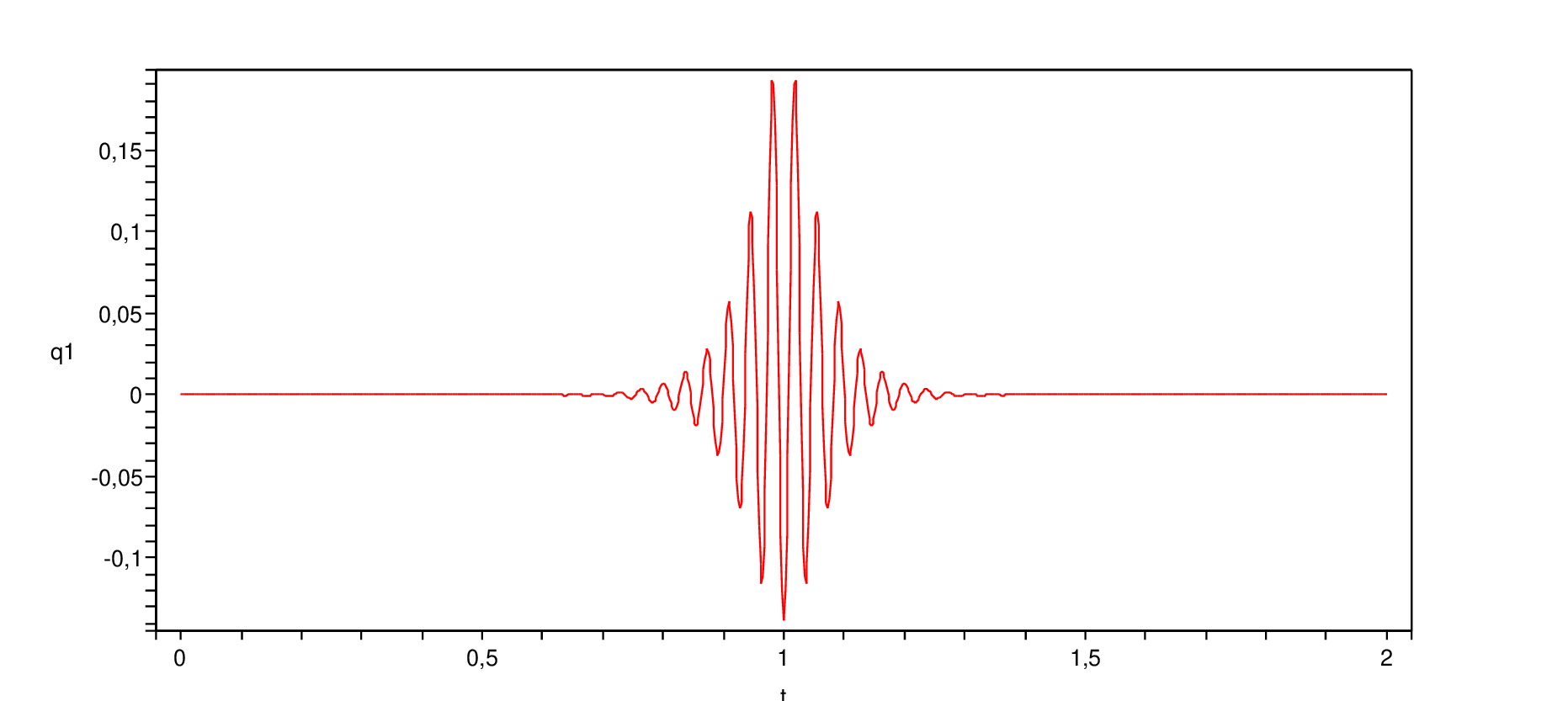}
  \end{center}
  \caption{\small Two primary symmetric homoclinic solutions 
of the scalar stationary GSHE ($\epsilon=-0.05$ and $\kappa=2$).} 
  \label{4:homoclinicorbits1}
\end{figure}

In the following we consider $u$ to be one dimensional 
and study stationary solutions of \eqref{4:GSHE} which satisfy the ordinary differential equation
\begin{equation}\label{4:SHE}
\epsilon u + \kappa u^2-u^3-(1+\partial^2_x)^2 u = 0\,.
\end{equation}
Obviously this equation has a reversible symmetry (if $u(x)$ satisfy the equation
then $u(-x)$ also does). It is well known that for small
negative $\epsilon$ this equation has two symmetric
homoclinic solutions~\cite{GL:95} similar to the ones shown on Figure~\ref{4:homoclinicorbits1}.
In this paper we study transversality of the homoclinic solutions,
which implies existence of multi-pulse homoclinic solutions and
a small scale chaos.

In order to describe the homoclinic phenomena it is convenient to rewrite
the equation \eqref{4:SHE} in the form of an equivalent Hamiltonian system~\cite{BelyakovGL1997,LLLB:04}:
 \begin{align}\label{4:HE}
\dot{q_1}&=q_2 & \dot{p_1}&=p_2-\epsilon q_1 - \kappa q_1^2 + q_1^3 \\
\notag \dot{q_2}&=p_2-q_1 & \dot{p_2}&=-p_1\,,
\end{align}
where the variables are defined by the following equalities
\begin{equation}\label{4:changeofvariables}
u=q_1,\quad u' = q_2,\quad -(u'+u''')=p_1 \quad \mathrm{and}\quad  u+u''=p_2
\end{equation}
and the Hamiltonian function has the form
\begin{equation}\label{4:H}
H_\epsilon=p_1 q_2 - p_2 q_1 + \frac{p_2^2}{2}+\epsilon \frac{q_1^2}{2}+\kappa \frac{q_1^3}{3}-\frac{q_1^4}{4}.
\end{equation}
The system \eqref{4:HE} is reversible with respect to the involution,
\begin{equation*}
S: (q_1,q_2,p_1,p_2) \rightarrow (q_1,-q_2,-p_1,p_2).
\end{equation*}
The origin is an equiblibrium of the system and the eigenvalues of the linearized vector field are
\begin{equation*}
\left\{\pm\sqrt{-1+\sqrt{\epsilon}},\,\pm\sqrt{-1-\sqrt{\epsilon}}\right\}\,.
\end{equation*} 
If $\epsilon<0$, the eigenvalues form a quadruple $\pm \beta_\epsilon\pm i\alpha_\epsilon$ where
\begin{eqnarray*}
\beta_\epsilon&=&\frac{\sqrt {2\sqrt {1-{\it \epsilon}}-2}}{2}=\sqrt{-\frac{\epsilon}4}\,(1+O(\epsilon))\,,\\
\alpha_\epsilon&=&\frac{\sqrt {2\sqrt {1-{\it \epsilon}}+2}}{2}=1+O(\epsilon)\,.
\end{eqnarray*}
At $\epsilon=0$ the eigenvalues collide forming two purely imaginary eigenvalues~$\pm i$ of multiplicity two. Moreover, 
the corresponding linearization of the vector field is not semisimple. 
Thus, the equilibrium point of system \eqref{4:HE} undergoes 
a Hamiltonian-Hopf bifurcation described in the book \cite{vdMeer1985} (see also \cite{Sok}).
In general position there are two possible scenarios of the bifurcation depending on the sign of a certain 
coefficient of a normal form. In the Swift-Hohenberg equation both scenarios are possible
and depend on the value of the parameter $\kappa$. In this paper we will consider the case when the equilibrium is stable at
the moment of the bifurcation (see \cite{Tre,LM} for more details) 
which corresponds to $|\kappa|>\sqrt{\frac{27}{38}}$ as shown in~\cite{BelyakovGL1997}. Also note that the degenerate case $|\kappa|=\sqrt{\frac{27}{38}}$ leads to some interesting phenomena including 
``homoclinic snaking" \cite{WoodC99,KnoblochW2008,CK:09}.

When $\epsilon<0$ is small, the equilibrium is a saddle-focus and the Stable Manifold Theorem 
implies the existence of two-dimensional stable $\mathbf{W}_\epsilon^{s}$ and unstable $\mathbf{W}_\epsilon^{u}$ 
manifolds for the equilibrium point. These manifolds are 
contained inside the zero energy level of the Hamiltonian $H_\epsilon$.

The original Hamiltonian \eqref{4:H} can be seen as a perturbation of an integrable Hamiltonian 
which can be derived from the normal form theory (see section \ref{Se:nf} for details). Since the normal form is integrable, 
its stable and unstable manifolds coincide
(see also discussion in \cite{GIMP:93} for the reversible set up).
In \cite{GL:95}, Glebsky and Lerman used the implicit function theorem
to prove the existence of two reversible (symmetric) homoclinic orbits for the original 
system \eqref{4:HE} when $\epsilon<0$ is small. 
As a matter of fact, this result follows from a more general 
study concerning a 1:1 resonance in four dimensional reversible vector fields 
(see \cite{GIMP:93}). Also the paper \cite{GL:95} conjectures that
the stable and unstable manifolds should intersect transversely yielding, 
in particular, the existence of countably many reversible homoclinic orbits. 
These orbits are known as \textit{multisolitons} for the Swift-Hohenberg equation 
and they have been the subject of study in several works (see \cite{ARC:98} and the references therein).

Note that no conclusion about the transversality of stable and unstable manifolds can be made using only the normal form theory. 
In this paper we study the splitting of the stable and unstable 
manifolds which happens beyond all orders of the normal
form theory. Let $\mathbf{p}_\epsilon$ be a symmetric homoclinic point 
belonging to one of the two primary symmetric homoclinic orbits.
In section \ref{Se:hi} we propose a natural way to select vectors $v^{u,s}_\epsilon$ tangent to
$\mathbf{W}_\epsilon^{s}$ and $\mathbf{W}_\epsilon^{u}$  
at $\mathbf{p}_\epsilon$ (see equation (\ref{Eq:vuvs})). The main goal of this paper
is to establish the following asymptotic formula for the value of the standard symplectic 
form on this pair of vectors:
\begin{equation}\label{4:Asymptoticformulahomoclinic}
\Omega(v^u_\epsilon,v^s_\epsilon)=e^{-\frac{\pi\alpha_\epsilon}{2\beta_\epsilon}}
\left(\omega_0(\kappa)+O(\epsilon)\right)\,.
\end{equation}
Note that $\mathbf{W}_\epsilon^{s}$ and $\mathbf{W}_\epsilon^{u}$ are two dimensional Lagrangian manifolds confined inside
the three dimensional energy level $\{H_\epsilon=0\}$.
These manifolds intersect along homoclinic orbits. Their intersection along the orbit of $\mathbf{p}_\epsilon$ is transverse 
(inside the energy level) provided $\Omega(v^u_\epsilon,v^s_\epsilon)\neq 0$. If $\omega_0(\kappa)\ne0$,
the  asymptotic formula (\ref{4:Asymptoticformulahomoclinic}) implies the transversality of the homoclinic orbit for small negative $\epsilon$,
and therefore $\omega_0(\kappa)$ is known as the \textit{splitting coefficient}.

We stress that the derivation of formula (\ref{4:Asymptoticformulahomoclinic}) does not 
rely substantially on the specific form of the Swift-Hohenberg equation and exactly the same asymptotic
expression (only the splitting coefficient may take different values) 
can be deduced {\em for a generic analytic family of reversible Hamiltonian systems undergoing a subcritical Hamiltonian-Hopf
bifurcation}. The details can be found in \cite{JP10}, where a majority of the arguments presented in this paper have been transformed into rigorous mathematical proofs. 

In this paper, the derivation of formula \eqref{4:Asymptoticformulahomoclinic} is not rigorous,
as it is based on numerous estimates and assumptions which are not proved here.
Nevertheless similar statements were proved in a similar context for
other problems \cite{G99,VGVL:01}.

\begin{figure}[t]
  \begin{center}
    \includegraphics[width=3.5in]{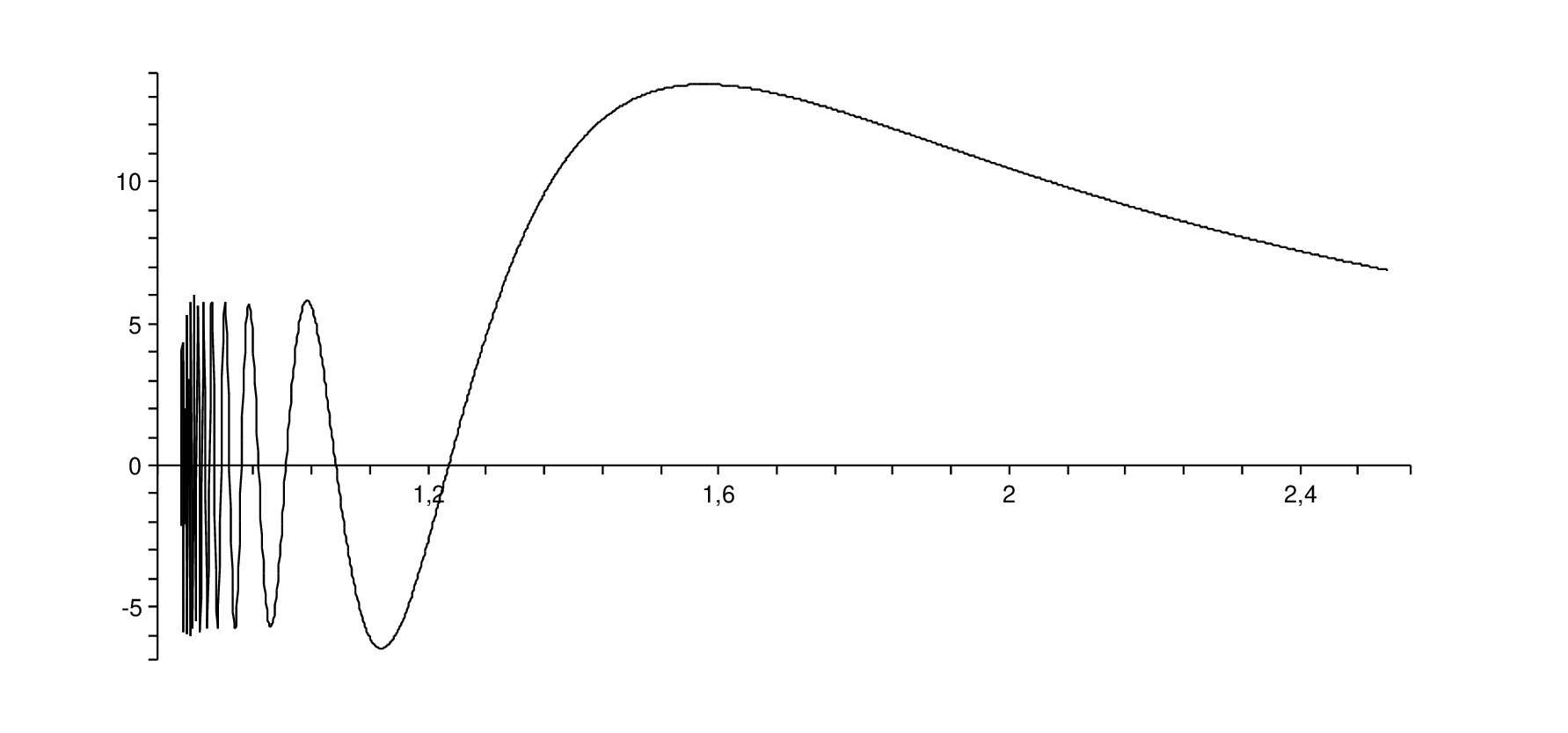}
  \end{center}
  \caption{Graph of the function $\mathrm{Im}(\Theta_0(\kappa))$ for $\kappa>\sqrt{\frac{27}{38}}$, $\omega_0(\kappa)=2\mathrm{Im}(\Theta_0(\kappa))$.} 
  \label{4:figStokesbeta}
\end{figure}

Supporting the validity of formula \eqref{4:Asymptoticformulahomoclinic} we perform a set of numerical experiments and compute the splitting coefficient using two distinct methods. This constant is related to a purely imaginary Stokes constant, and
Figure~\ref{4:figStokesbeta} gives an idea about its behaviour
as a function of the parameter $\kappa$. The value of the Stokes constant comes from the study of the Hamiltonian \eqref{4:H}
at the exact moment of the bifurcation (i.e. at $\epsilon=0$). We will discuss the relevant definitions in section~\ref{Se:Stokes}
and some methods for its numerical evaluation in section~\ref{Se:num}.

Recently, S. J. Chapman and G. Kozyreff \cite{CK:09}  used the multiple-scales analysis beyond all orders
to study localised patterns  which emerge from a subcritical modulation instability in the Swift-Hohenberg equation. 
Their analysis captured exponentially small phenomena by means of optimal truncation of certain formal expansions 
combined with a study of their analytical continuation in a vicinity of the Sokes lines. 
Technically our approach is different and we do not require higher order terms, additionally our
approach has the advantage of being directly applicable to study the exponentially small splitting of invariant manifolds 
near generic Hamiltonian-Hopf bifurcations, for which the Swift-Hohenberg is a particular example. 

The rest of the paper is organised in the following way. In Section~\ref{Se:hi}
we discuss the definition of an homoclinic invariant which provides a very convenient
tool for measuring the splitting of invariant manifolds. In Section~\ref{Se:nf}
we review some facts from the normal form theory which will be necessary
for the exposition of our results. Section~\ref{Se:Stokes} contains the definition of the Stokes constant. An informal derivation of the asymptotic formula \eqref{4:Asymptoticformulahomoclinic}
which describes the splitting of invariant manifolds of the
stationary Swift-Hohenberg equation near the Hamiltonian-Hopf bifurcation is placed in Section~\ref{Se:asymp}.
%In particular this formula implies existence of two transversal 
%primary homoclinic orbits and gives a sharp lower (exponentially small) bound for the splitting of the invariant manifolds.
As the derivation of the asymptotic formula is not rigorous
we perform a set of high-precision numerical experiments
in order to confirm its validity. Moreover, similar to
many other problems which involve exponentially small splitting of
invariant manifolds \cite{VGVL:01}, the asymptotic formula contains a splitting 
coefficient which comes from an auxiliary problem and requires numerical evaluation. 
The results of our numerical experiments are reported in Sections~\ref{Se:num} and~\ref{Se:num3}.

\subsection{Homoclinic invariant\label{Se:hi}}
In a study of homoclinic trajectories, both numerical and analytical, 
it is usually important to have a convenient basis in the tangent
space to the stable and unstable manifolds. Below we provide a definition
adapted to our problem. This definition can be of independent interest 
as it can be easily extended onto hyperbolic equilibria of higher
dimensional systems (not necessarily Hamiltonian).

Suppose that the origin is an equilibrium of a Hamiltonian vector field $X_H$
and that $\pm\beta\pm i\alpha$ are the eigenvalues of $DX_H(0)$. 
Then the origin has a two dimensional stable manifold. According to Hartman~\cite{MR0141856} the restriction of 
the vector field on $W^s_{loc}$ can be linearised by a $C^1$ change of variables.
In the polar coordinates the linearised dynamics on $W^s_{loc}$ takes the form:
$$
\dot r=-\beta r\,\qquad \dot\varphi=\alpha\,.
$$
It is convenient to introduce $z=-\ln r$ so that
$$
\dot z=\beta\,.
$$ 
Then the local stable manifold is the image of a function 
$$
\Gamma^s:\{(\varphi,z):\varphi\in S^1,z>-\log r_0\}\to\mathbb R^4
$$
where $r_0$ is the radius of the linearisation domain and $S^1$ is the unit circle.
Since $\Gamma^s$ maps trajectories into trajectories we can propagate it uniquely
along the trajectories of the Hamiltonian system using the property
\begin{equation}\label{Eq:Gammau}
\Gamma^s(\varphi+\alpha t,z+\beta t)=\Phi^t_H\circ \Gamma^s(\varphi,z)
\end{equation}
where $\Phi^t_H$ is the flow defined by the Hamiltonian equation.
Note that 
$$
\Gamma^s(\varphi+2\pi,z)=\Gamma^s(\varphi,z)
$$
since $\varphi$ is the angle component of the polar coordinates.
Moreover,
$$
\lim_{z\to+\infty}\Gamma^s(\varphi,z)=0\,.
$$
Differentiating $\Gamma^s$ along a trajectory we see that it
satisfies the non-linear PDE:
\begin{equation}\label{Eq:mainPDE}
\alpha\partial_\varphi\Gamma+\beta\partial_z\Gamma=X_H(\Gamma)\,.
\end{equation}
Each of the derivatives $\partial_z\Gamma^s$ and $\partial_\varphi\Gamma^s$ defines
a vector field on $W^s$. 
%It can be checked that these vector fields are defined uniquely.
%Indeed, since $\Gamma^s$ is defined uniquely up to a translation of $(z,\varphi)$ plane,
%the vector fields are independent from the freedom in the definition of $\Gamma^s$.
The equation (\ref{Eq:Gammau}) implies that $\partial_z\Gamma^s$ and $\partial_\varphi\Gamma^s$  
are invariant under the restriction of the flow $\Phi^t_H\Bigr|_{W^s}$. 

We can define $\Gamma^u$ applying the same arguments to the Hamiltonian $-H$. In this case it is convenient
to set $z=\ln r$ to ensure that $\Gamma^u$ satisfies the same PDE as $\Gamma^s$.
In a reversible system with a reversing involution $S$, it is convenient to
set 
\begin{equation}\label{Eq:unstablemfd}
\Gamma^u(\varphi,z)=S\circ\Gamma^s(-\varphi,-z).
\end{equation}
Now suppose that the system has a homoclinic trajectory $\gamma_h$. Let us choose a point $\mathbf{p}_h\in\gamma_h$.
The freedom in the definition allows us
to assume that $\mathbf{p}_h=\Gamma^s(0,0)=\Gamma^u(0,0)$ without loosing in generality.
This condition completely eliminates the freedom from the definition of $\Gamma^u$ and $\Gamma^s$.

In a Hamiltonian system the symplectic form provides a natural tool for studying transversality
of invariant manifolds. Thus we arrive at the following,
\begin{definition}[Homoclinic Invariant]The homoclinic invariant $\omega$ is defined by the formula,
\begin{equation}\label{Eq:hominv}
\omega=\Omega(\partial_\varphi\Gamma^u(0,0),\partial_\varphi\Gamma^s(0,0))\,.
\end{equation}
\end{definition}
This definition is a natural extension of the homoclinic invariant defined for homoclinic
orbits of area-preserving maps~\cite{VGVL:01}.

In the left hand side of the asymptotic formula (\ref{4:Asymptoticformulahomoclinic}) 
we use the notation 
\begin{equation}\label{Eq:vuvs}
v_\epsilon^{u,s}=\partial_\varphi\Gamma^{u,s}(0,0).
\end{equation}
It is easy to see that $\omega$ takes the same 
value for all points of the homoclinic trajectory
$\gamma_h=\{\Phi^t_H(\mathbf{p}_h):t\in\mathbb R\}$. Indeed it follows from \eqref{Eq:Gammau} that
$$
\partial_\varphi\Gamma^s(\alpha t,\beta t)=D\Phi^t_H(\mathbf{p}_\epsilon)\partial_\varphi \Gamma^s(0,0),
$$
and a similar identity is valid for the unstable manifold. Since the Hamiltonian flow $\Phi^t_H$ is symplectic, 
we conclude that 
$\Omega(\partial_\varphi\Gamma^u(\alpha t,\beta t),\partial_\varphi\Gamma^s(\alpha t,\beta t))
=\Omega(\partial_\varphi\Gamma^u(0,0),\partial_\varphi\Gamma^s(0,0))=\omega$.

Since $\Gamma^s$ and $\Gamma^u$ are Lagrangian and belong to the energy level $H=H(0)$, which is three-dimensional,
the inequality $\omega\ne0$ implies the transversality of the homoclinic trajectory. 
Indeed, if $\omega\ne0$, the vectors $\partial_\varphi\Gamma^u(0,0)$, $\partial_\varphi\Gamma^s(0,0)$
and $X_H(\mathbf{p}_h)$ are linearly independent and therefore span the tangent space to
the energy level at $\mathbf{p}_h$.

We note that we can define two vectors tangent to $W^s$ and
another two vectors tangent to $W^u$ at $\mathbf{p}_h\in W^s\cap W^u$.
So we could use 
\begin{equation}\label{Eq:xy}
\omega_{x,y}:=\Omega(\partial_x\Gamma^u(0,0),\partial_y\Gamma^s(0,0)),\quad x,y\in\left\{\varphi,z\right\}
\end{equation}
instead of $\omega$. But these invariants are not independent. Indeed, 
$$
\alpha\partial_\varphi\Gamma^u(0,0)+\beta\partial_z\Gamma^u(0,0)=\alpha\partial_\varphi\Gamma^s(0,0)+\beta\partial_z\Gamma^s(0,0)
$$
as both expressions are equal to $X_H(\mathbf{p}_h)$. Then 
equation (\ref{Eq:xy}) implies
$$
\alpha^2\omega-\beta^2\omega_{z,z}=0,\quad \alpha\omega+\beta\omega_{\varphi,z}=0,\quad \alpha\omega+\beta\omega_{z,\varphi}=0.
$$
In the derivation of these identities it is also necessary to take into account that
$W^{u,s}$ are Lagrangian (i.e., the symplectic form $\Omega$ vanishes on their tangent spaces).

In the case of the Swift-Hohenberg equation the system of PDE (\ref{Eq:mainPDE})
can be conveniently replaced by a single scalar PDE of higher order
obtained from (\ref{4:SHE}) by replacing $\partial_x$ with the differential operator
$$
\partial=\alpha_\epsilon\partial_\varphi + \beta_\epsilon\partial_z.
$$
Let us use $u^{\pm}_\epsilon$ to denote the first component of $\mathbf{\Gamma}^u_\epsilon$
and $\mathbf{\Gamma}^s_\epsilon$ respectively, then $u^\pm_\epsilon$ satisfies the equation
\begin{equation}\label{Eq:SHPDE}
(1+\partial^2)^2 u = \epsilon u + \kappa u^2-u^3\,.
\end{equation}
Its other components can be restored using (\ref{4:changeofvariables}). 
The Swift-Hohenberg equation is reversible and following \eqref{Eq:unstablemfd} we define
$$
u^{+}_\epsilon(\varphi,z)=u^{-}_\epsilon(-\varphi,-z)\,.
$$
We also assume that $\mathbf{\Gamma}^s_\epsilon(0,0)=\mathbf{\Gamma}^u_\epsilon(0,0)$ is 
the primary symmetric homoclinic point. Then the formula
for the homoclinic invariant can be rewritten in terms of $u^{-}$:
%$$
%
%\partial u^{+}=-\partial u^{-} , \partial^2 u^{+}=\partial^2 u^{-} and \partial^3 u^{+}=-\partial^3 u^{-}
%
%Then,
\begin{equation}
\omega=2\partial_\varphi\left((u^{-})^2 +u^{-}\partial^2u^{-})\right)
\end{equation}
where the derivatives are  evaluated at $(\varphi,z)=(0,0)$.

\subsection{Normal form of the Swift-Hohenberg equation\label{Se:nf}} 
The most convenient description of the bifurcation is obtained with the help of the normal form.
As a first step the quadratic part of the Hamiltonian \eqref{4:H}
is normalised with the help of a linear symplectic transformation (similar to  \cite{NBRC:74}): 
\begin{equation*}
T= \left( \begin {array}{cccc} 0&-1/4\,\sqrt {2}&-1/2\,\sqrt {2}&0
\\\noalign{\medskip}1/4\,\sqrt {2}&0&0&1/2\,\sqrt {2}
\\\noalign{\medskip}\sqrt {2}&0&0&0\\\noalign{\medskip}0&-\sqrt {2}&0&0
\end {array} \right)
\end{equation*}
which transforms \eqref{4:H} into
\begin{equation}\label{4:H2}
\begin{split}
H_\epsilon=&-(q_2p_1-q_1p_2)+\frac{1}{2}(q_1^2+q_2^2)+\frac{1}{4}p_1^2\epsilon
-\frac{\sqrt{2}}{12}\kappa p_1^3+\frac{1}{4}q_2p_1\epsilon-\frac{\sqrt{2}}{8}\kappa q_2p_1^2+\\
&\frac{1}{16}q_2^2\epsilon-\frac{\sqrt{2}}{16}\kappa q_2^2p_1
-\frac{\sqrt{2}}{96}\kappa q_2^3-\frac{1}{16}p_1^4-\frac{1}{8}q_2p_1^3
-\frac{3}{32}q_2^2p_1^2-\frac{1}{32}q_2^3p_1-\frac{1}{256}q_2^4
\end{split}
\end{equation}
where we keep the same notation for the variables.
Note that the involution $S$ in the new coordinates takes the form
\begin{equation}\label{4:Involution}
\tilde{S}: (q_1,q_2,p_1,p_2) \rightarrow (-q_1,q_2,p_1,-p_2).
\end{equation}
Now, with the quadratic part in normal form, we can apply 
the standard normal form procedure to normalize the Hamiltonian \eqref{4:H2} 
up to any order: There is a near identity canonical change of variables $\Psi_n$ 
which normalizes all terms of order less than equal to $n$
and transforms the Hamiltonian to the following form:
%\begin{equation}\label{4:HN}\begin{split}
%H_\epsilon=I_1 - \frac{1}{2}I_2-\frac{1}{8}\epsilon I_1 - \frac{1}{8}\epsilon I_3 
%+ \left( {\frac {7}{1296}}\,{\kappa}^{2}+\frac{1}{32} \right) %{I_1}^{2}+ \left( {\frac {3}{64}}-{\frac {65}{
%864}}\,{\kappa}^{2} \right) I_1{I_3}\\+ \left( -{\frac {19}{576}}\,
%{\kappa}^{2}+{\frac {3}{128}} \right) {I_3}^{2} + F
%\end{split}
%\end{equation}

\begin{equation}\label{4:HN}
H_\epsilon=H_\epsilon^n + \mathrm{higher}\ \mathrm{order}\ \mathrm{terms}
\end{equation}
where
\begin{equation*}
H_\epsilon^n=-I_1 + I_2+\sum_{\substack{3i+2j+2l\geq4\\i+j\geq1}}^n a_{i,j,l}I_1^iI_3^j\epsilon^l
\end{equation*}
with 
\begin{equation*}
I_1=q_2p_1-q_1p_2,\qquad I_2=\frac{q_1^2+q_2^2}{2},\qquad I_3=\frac{p_1^2+p_2^2}{2}.
\end{equation*}
This normalization preserves the reversibility with respect to the involution \eqref{4:Involution}. 
In the case of the GSHE the normal form up to the order five has the form (see Appendix~\ref{Ap:A}
for more details about the change of variables)
\begin{equation*}
H_\epsilon^5=-I_1 + \left(I_2+\frac{1}{4}\epsilon I_3 + \eta I_3^2\right) 
+\left(\frac{1}{8}\epsilon I_1+\mu\, I_1{I_3}\right)\,.
\end{equation*}
The leading part of the normal form includes two parameters which can be explicitly
expressed in terms of the original parameter $\kappa$:
\begin{equation*}
\eta=4\left( {\frac {19}{576}}{\kappa}^{2}-{\frac {3}{128}} \right)\quad \mathrm{and}\quad 
\mu=2\left( {\frac {65}{864}}\,{\kappa}^{2} -{\frac {3}{64}}\right)\,.
\end{equation*}
The geometry of the invariant manifolds depends on the sign of $\eta$ \cite{vdMeer1985}.
In the case of GSHE, if 
\begin{equation*}
\left|\kappa\right|>\sqrt{\frac{27}{38}},
\end{equation*}
then $\eta>0$ \cite{GL:95}, and the truncated normal form has a continuum of homoclinic orbits
among which exactly two are reversible, i.e., symmetric with respect to the
involution \eqref{4:Involution}. 

In order to describe the geometry of the invariant manifolds near the bifurcation
it is convenient to introduce the new parameter $\epsilon=-4\delta^2$ and 
perform the standard scaling:
\begin{equation*}
q_1=\delta^2Q_1,\quad q_2=\delta^2Q_2,\quad p_1=\delta P_1,\quad p_2=\delta P_2\,.
\end{equation*}
This change of variables is not symplectic, nevertheless  it preserves the form of the Hamiltonian equations
since the symplectic form gains a constant factor $\delta^3$, so we have to multiply the Hamiltonian
by $\delta^{-3}$ in order to return back to the standard symplectic form.
 The Hamiltonian $H_\epsilon^n$ is transformed into,
\begin{equation*}
h_\delta^n=- \mathcal{I}_1 + \left(\mathcal{I}_2-\mathcal{I}_3+\eta\mathcal{I}_3^2\right)\delta
+ \left(-\frac{1}{2} \mathcal{I}_1+\mu\, \mathcal{I}_1{\mathcal{I}_3}\right)\delta^2+O(\delta^3),
\end{equation*}
where the $\mathcal{I}_i$'s are defined in the same way as the $I_i$'s but in the new variables $Q$ and $P$. 
This Hamiltonian system has an equilibrium at the origin characterized by a quadruple of complex eigenvalues
$\pm i\alpha_{n,\epsilon}\pm \beta_{n,\epsilon}$, where 
$\alpha_{n,\epsilon}=1+\frac{1}{2}\delta^2+O(\delta^4)$ and $\beta_{n,\epsilon}=\delta-\frac{1}{2}\delta^3+O(\delta^5)$.

The equilibrium has a two dimensional stable and two dimensional unstable manifolds. Thus, following \eqref{Eq:mainPDE} we parametrize these manifolds by solutions
of the partial differential equation:
\begin{equation}\label{4:HEN}
\left(\alpha_{n,\epsilon}\partial_\varphi + \beta_{n,\epsilon}\partial_z\right)\mathbf{\Upsilon}_n
=X_{h_\delta^n}(\mathbf{\Upsilon}_n).
\end{equation}
The function $\mathbf{\Upsilon}_n(\varphi,z)$ is real-analytic, converges to zero as $z\to\pm\infty$ 
and is $2\pi$-periodic in $\varphi$. 
Taking into account the rotational symmetry of the normal form Hamiltonian, we can look for the
solution of this equation in the form:
\begin{equation*}
\begin{split}
\mathbf{\Upsilon}_n(\varphi,z)&=\bigl(R_n(z)\cos(\theta_n(\varphi,z)),R_n(z) \sin(\theta_n(\varphi,z)),\\
&\qquad r_n(z)\cos(\theta_n(\varphi,z)),r_n(z)\sin(\theta_n(\varphi,z))\bigr)
\end{split}
\end{equation*}
\normalsize
where $R_n(z)$, $r_n(z)$ and $\theta_n(\varphi,z)$ are real analytic functions.
In particular, for $n=5$ it is not difficult to see that the eigenvalues of $DX_{h_\delta^5}(0)$ are the quadruple $\pm\beta_{5,\epsilon}\pm i\alpha_{5,\epsilon}$ where,
$$
\beta_{5,\epsilon}=\delta\qquad \alpha_{5,\epsilon}=1+\frac{\delta^2}2\,.
$$
Thus, we get the following system of equations:
\begin{equation*}
\begin{split}
\beta_{5,\epsilon}R_5'=-\delta r_5\left(1-\eta r_5^2\right)
\,,
\qquad
\beta_{5,\epsilon}r_5'=-\delta R_5\,,\qquad
\\
\left(\alpha_{5,\epsilon}\partial_\varphi + \beta_{5,\epsilon}\partial_z\right)
\theta_5=
1+\frac{\delta^2}2(1-\mu r_5^2)\,.
\end{split}
\end{equation*}
From these equations we conclude that
\begin{equation*}
\begin{split}
r_5=\sqrt{\frac{2}{\eta}}\frac1{\cosh z}\,,\qquad
R_5=\sqrt{\frac{2}{\eta}}\frac{\sinh z}{\cosh^2 z},\\
\theta_5=\varphi-\frac{\delta^2\mu}2\int^zr_5^2dz
=\varphi-\frac{\delta\mu}{\eta}\frac{\sinh z}{\cosh z}
\,.
\end{split}
\end{equation*}
We see that $(r_5(z),R_5(z))$ runs over a homoclinic loop when $z$ varies
from $-\infty$ to~$+\infty$.

In general the parameterization $\mathbf{\Upsilon}_n$ is the unique solution of 
\eqref{4:HEN} such that $R_n(0)=0$ and $\theta_n(\varphi,0)=\varphi$. 
Thus, $\mathbf{\Upsilon}_n(\varphi,z)$ belongs to the symmetry plane associated 
with the involution \eqref{4:Involution} if and only if $z=0$ and $\varphi=0$ 
or $\varphi=\pi$. Therefore, there are exactly 2 symmetric homoclinic points. 
Let us call these homoclinic orbits the primary reversible homoclinic orbit.

\subsection{Stokes constant}\label{Se:Stokes}

In this subsection we define the Stokes constant for the GSHE at $\epsilon=0$. 
Although the equilibrium at the origin is not hyperbolic (its eigenvalues are $\pm i$ with multiplicity two), it still has
invariant manifolds \cite{JP10} which can be non-real.
More precisely, we look for complex analytic solutions of the following equation
\begin{equation}\label{Eq:inner}
(1+(\partial_\varphi + \partial_\tau)^2)^2 u = \kappa u^2-u^3\,,
\end{equation}
which decay polynomially in a sectorial neighbourhood of infinity in the $\tau$ variable
 and which are $2\pi$-periodic in $\varphi$.
These solutions parametrize a certain complex stable (unstable) invariant manifold of the origin which is immersed in $\mathbb{C}^4$. 
In \cite{JP10} it is shown (for similar problems see \cite{VGVL:01,Baldoma06,OliveSS03}) that equation \eqref{Eq:inner} has an analytic solution $u=u_0^-$ with the following asymptotic behaviour:
\begin{equation*}
u^{-}_0(\varphi,\tau)=\frac{P_1(\varphi)}{\tau} +\frac{P_2(\varphi)}{\tau^2}+O(\tau^{-3})
\end{equation*}
in the set 
$$\tau\in\mathcal{D}^{-}_{r,\theta_0}=\left\{\tau\::\:\left|\arg(\tau+r)\right|>\theta_0\right\}
\,,$$
where $\theta_0$ is a small fixed constant, $r$ is sufficiently large and
\begin{equation}\label{Eq:p1p2}
P_1=\frac{i\cos \left( \varphi \right)}{\sqrt {\eta}}\,,
\qquad
P_2=
\frac{i}{\sqrt {\eta}} \left( \frac{\mu}{{\eta}}+\frac{1}{2}\right)\sin(\varphi)-{\frac {\kappa\,\cos \left( 2\,\varphi \right) }{18\eta}}-{
\frac {\kappa}{2\eta}}\,.
\end{equation}
The function $u^-_0$ is $2\pi$-periodic in $\varphi$. More generally it is possible to prove (see \cite{JP10}) that there exist unique trigonometric polynomials $P_k$ for $k\geq 3$ of degree $k$ satisfying $P_k(\varphi)=(-1)^k\overline{P_k(-\overline\varphi)}$ such that $\hat{u}_0(\varphi,\tau):=\sum_{k\geq1}P_k(\varphi)\tau^{-k}$ solves formally equation \eqref{Eq:inner} and moreover,
\begin{equation*}
u^{-}_0(\varphi,\tau)=\sum_{k=1}^{N}P_k(\varphi)\tau^{-k}+O(\tau^{-(N+1)}).
\end{equation*}
Taking into account \eqref{Eq:p1p2} we have that $\hat{u}_0(\varphi,\tau)=\overline{\hat{u}_0(-\overline\varphi,-\overline\tau)}$ and the unique formal solution $\hat{u}_0$ is known as the \textit{formal separatrix}.

Equation \eqref{Eq:inner} has a second solution $u=u^+_0$ with
$$
u^+_0(\varphi,\tau)=\overline{u^-_0(-\overline\varphi,-\overline\tau)}\,.
$$
It has the same asymptotic behaviour as $u^-_0$ but is defined in a different sector, more precisely,
it is defined for $\tau$ such that $-\overline\tau\in\mathcal{D}^{-}_{r,\theta_0}$. The solutions $u^\pm_0$ have a
common asymptotics on the intersection of their domains but they do not
typically coincide. The difference of these two solutions can be described in the following way.
We can restore 4-dimensional vectors $\mathbf{\Gamma}_0^{\pm}$ using equations 
\eqref{4:changeofvariables} with $'$ replaced by $\partial_\varphi+\partial_\tau$.
In particular, the first component of $\mathbf{\Gamma}_0^{\pm}$ coincides with $u^\pm_0$. 
The functions $\mathbf{\Gamma}^{\pm}_0$ are parameterizations of the stable and unstable manifolds 
and satisfy the following non-linear partial differential equation,
\begin{equation}\label{4:VFH0}
(\partial_\varphi+\partial_\tau)\mathbf{\Gamma}^{\pm}_0=X_{H_0}(\mathbf{\Gamma}^{\pm}_0),
\end{equation}
where $H_0$ denotes the Hamiltonian \eqref{4:H} at the exact moment of bifurcation $\epsilon=0$.
Let
$$
\Delta_0(\varphi,\tau)=
\mathbf{\Gamma}_0^{+}(\varphi,\tau)-\mathbf{\Gamma}_0^{-}(\varphi,\tau)
$$
and
$$
\theta_0(\varphi,\tau)=
\Omega\bigl(\Delta_0(\varphi,\tau),\partial_\varphi\mathbf{\Gamma}_0^{+}(\varphi,\tau)\bigr)\,,
$$
where $\Omega$ is the standard symplectic form. In \cite{JP10} it is proved that there is a constant $\Theta_0(\kappa)$ such that
\begin{equation}\label{Eq:theta0}
\theta_0(\varphi,\tau)
=\Theta_0(\kappa)e^{-i(\tau-\varphi)}+O(e^{-(2-\epsilon_0)i(\tau-\varphi)})
\end{equation} 
as $\mathrm{Im}\,\tau \rightarrow-\infty$ and for very small $\epsilon_0>0$.
The constant $\Theta_0(\kappa)$ is known as the\/ {\em Stokes constant}.
The Stokes constant can be defined by the following limit:
\begin{equation}\label{4:SC}
\Theta_0(\kappa):=\lim_{\mathrm{Im}(\tau)\rightarrow -\infty} \theta_0(\varphi,\tau)e^{i(\tau-\varphi)}\,.
\end{equation}
We note that the value of the Stokes constant cannot be obtained 
from our arguments. Fortunately the numerical evaluation of this constant is reasonably easy. 
Figure~\ref{4:figStokesbeta} shows the values of $\mathrm{Im}\,\Theta_0(\kappa)$ plotted
against $\kappa$ for $\kappa>\kappa_0=\sqrt{\frac{27}{38}}$. The picture suggests that the Stokes constant vanishes 
infinitely many times and that its zeros accumulate to $\kappa_0$.

\subsection{Asymptotic formula for the homoclinic invariant\label{Se:asymp}}

In this section we derive the asymptotic formula (\ref{4:Asymptoticformulahomoclinic}) for the homoclinic invariant
of the primary symmetric homoclinic orbit. Our method is not rigorous
and relies on the complex matching approach
similar to one used for the standard map and the rapidly perturbed pendulum (see \cite{VGVL:01}).
We point out that in the latter two cases the method leaded to a complete proof of
asymptotic formulae similar to (\ref{4:Asymptoticformulahomoclinic}). Our approach has certain 
similarity to the complex matching methods used in \cite{HakimM93,CK:09} but is different in several
important technical details.

At the end of section \ref{Se:nf} we obtained an approximation of the separatrix in the normal form coordinates.
Transforming $\Upsilon_5(\varphi,z)$ back to the original coordinates we obtain 
the following approximation:
%
%Let us formally derive an asymptotic formula for the homoclinic invariant 
%of one of the two primary homoclinic orbits. Denote by $\pi_i:\mathbb{C}^4\rightarrow\mathbb{C}$ be the standard projection onto the $i$ component. There exists 
%a parameterisation $\mathbf{\Gamma}^u_\epsilon(\varphi,z)$ of the unstable manifold 
%$\mathbf{W}_\epsilon^{u}$ satisfying the following PDE,
%\begin{equation*}
%(\alpha_\epsilon\partial_\varphi + \beta_\epsilon\partial_z)\mathbf{\Gamma}^u_\epsilon= X_{H_\epsilon}(\mathbf{\Gamma}^u_\epsilon)
%\end{equation*}
%and the boundary conditions,
%\begin{equation*}
%\mathbf{\Gamma}^u_\epsilon(\varphi+2\pi,z)=\mathbf{\Gamma}^u_\epsilon(\varphi,z),\  \ 
%\lim_{z\rightarrow -\infty} \mathbf{\Gamma}^u_\epsilon(\cdot,z)=0,\ \ \mathbf{\Gamma}^u_\epsilon(0,0) \in \mathrm{Fix}(S)
%\end{equation*}
%Moreover, the first component of $\mathbf{\Gamma}^u_\epsilon$ has the following asymptotic behavior,
%
\begin{eqnarray}\label{4:leadingorderGamma_epsilon}
\lefteqn{u^{-}_\epsilon(\varphi,z)=
-\frac{1}{\sqrt {\eta}}{\frac {\cos \left( \varphi \right)}{\cosh
 \left( z \right) }}\delta
}
\\
\nonumber
&&\quad+\left(\frac{9\kappa+\kappa\cos(2\varphi)}{18\eta}\frac{1}{\cosh^2(z)}
-\frac{1}{\sqrt{\eta}}\left(\frac{\mu}{\eta}+\frac{1}{2}\right)\frac{\sin(\varphi)\sinh(z)}{\cosh^2(z)}\right)\delta^2+O(\delta^3)
\end{eqnarray}
where $\epsilon=-4\delta^2$. Since the function in the right-hand-side of the equation is even, 
it also approximates the stable separatrix represented by $u^{+}_\epsilon(\varphi,z)=u^{-}_{\epsilon}(-\varphi,-z)$.
A more accurate approximation with a $O(\delta^n)$ error can be obtained with the help of higher order normal form
theory, but naturally none of those approximations can distinguish between the stable and unstable separatrices
and we come to the conclusion that $$u^{-}_\epsilon(\varphi,z)-u^{+}_\epsilon(\varphi,z)=O(\delta^n)$$
for all $n$. Of course the constant in this upper bound may depend on the point $(\varphi,z)$.
%
%
%Due to \eqref{4:leadingorderGamma_epsilon} we know that,
%\begin{equation*}
%\pi_1\left(\mathbf{\Gamma}^u_\epsilon(\varphi,z)-\mathbf{\Gamma}^s_\epsilon(\varphi,z)\right)=O(\delta^{3})
%\end{equation*}
%and in fact this different is smaller than any fixed power of $\delta$ and valid for $(\varphi,z)\in\mathbb{R}\times(-L,L)$ for some $L>0$. 
%
Therefore, the difference between the stable and unstable parametrisation cannot be detected using power series of the perturbation theory,
and we say it is beyond all algebraic orders. 
A rather standard approach to the problem is based on studying the analytical continuation 
of the parametrisations and looking for places in the complexified variables 
where the leading orders of the approximation \eqref{4:leadingorderGamma_epsilon} grow significantly. 
We note that the variables $z$ and $\varphi$ play different roles,
in particular we assume that $\varphi$ is kept real or, more precisely,
in a fixed narrow strip around the real axis.

It is easy to see that the leading orders of $u^-_\epsilon$ have poles 
at $z=i\frac{\pi}{2}+ki\pi$ for any integer $k$. In the following we study the behaviour 
of the parametrisations near the singular point $z=i\frac{\pi}{2}$. The first step is to re-expand 
the functions in Laurent series around the singularity and introduce a new variable
\begin{equation}\label{4:ztaurelation}
\tau=\frac{\alpha_\epsilon}{\beta_\epsilon}z-i\frac{\pi\alpha_\epsilon}{2\beta_\epsilon}.
\end{equation}
Substituting this new variable into \eqref{4:leadingorderGamma_epsilon} 
and expanding around $\tau=0$ we conclude that 
\begin{equation}\label{4:leadingorderGamma_0}
u^{-}_\epsilon(\varphi,\tfrac{\beta_\epsilon}{\alpha_\epsilon}\tau+i\tfrac{\pi}{2})
=\left(\frac{P_1(\varphi)}{\tau} 
+\frac{P_2(\varphi)}{\tau^2}+O(\tau^{-3})\right)+O(\epsilon)
\end{equation}
where $P_1$ and $P_2$ are the same as in \eqref{Eq:p1p2}
and the error terms come from the analysis of the next order corrections.
In this analysis we consider the terms in \eqref{4:leadingorderGamma_epsilon} 
which are most divergent and in this way obtain the essential behaviour of $u^{-}_\epsilon$ 
around the singularity. 

Transforming the equation \eqref{Eq:SHPDE} to the variable \eqref{4:ztaurelation}, setting $\epsilon=0$
and noting that $\alpha_0=1$, we obtain equation \eqref{Eq:inner} considered in the previous subsection.
%
%Using this information around the singularity we can select an unique analytic function,
%$$
%\mathbf{\Gamma}_0^{-}:\mathcal{S}_h\times\mathcal{D}^{-}_{r,\theta_0}\rightarrow\mathbb{C}^4
%$$
%for some $0<\theta_0<\frac{\pi}{4}$ and $h,r>0$ which is $2\pi$-periodic in $\varphi$ and is a solution of the equation,
%\begin{equation}
%\partial_\varphi\mathbf{\Gamma}_0^{-}+\partial_\tau\mathbf{\Gamma}_0^{-}=X_{H_0}(\mathbf{\Gamma}_0^{-})
%\end{equation}
%such that,
%$$
%\pi_1\left(\mathbf{\Gamma}_0^{-}(\varphi,\tau)\right)=\frac{i}{\sqrt {\eta}}\cos \left( \phi \right)\frac{1}{\tau} +\left(i\frac{1}{\sqrt {\eta}} \left( \frac{\mu}%%{{\eta}}+\frac{1}{2}\right)\sin(\phi)-{\frac {\kappa\,\cos \left( 2\,\phi \right) }{18\eta}}-{
%\frac {\kappa}{2\eta}}\right)\frac{1}{\tau^2}+\mathcal{O}(\tau^{-3})
%$$
%
The following method is known as ``complex matching" and is based on the observation
that $u^\pm_0$ approximate $u^\pm_\varepsilon$ in a region where $|z-i\frac\pi2|$ is small
but $\tau$ is still large. 
Taking into account \eqref{4:leadingorderGamma_0} we conclude that
\begin{eqnarray}\label{Eqs:upmapp}
u^-_\epsilon(\varphi,\tfrac{\beta_\epsilon}{\alpha_\epsilon}\tau+i\tfrac{\pi}{2})&=&u_0^{-}(\varphi,\tau)+O(\epsilon)\,,\\
u^+_\epsilon(\varphi,\tfrac{\beta_\epsilon}{\alpha_\epsilon}\tau+i\tfrac{\pi}{2})&=&u_0^{+}(\varphi,\tau)+O(\epsilon)\,.
\end{eqnarray}
in a neighbourhood of a segment of the imaginary axis where $\Im\tau$ is large negative.
In a rigorous justification of the method we use the interval $-R\log\epsilon^{-1}<\Im\tau<-R$, where $R$ is a large constant.

Now restoring the 4-dimensional vectors $\mathbf{\Gamma}^{u,s}_\epsilon$ using the relations \eqref{4:changeofvariables} we obtain the following estimate for the difference,
\begin{equation}\label{Eq:est1}
\Delta(\varphi,\tfrac{\beta_\epsilon}{\alpha_\epsilon}\tau+i\tfrac{\pi}{2})=-\Delta_0(\varphi,\tau)+O(\epsilon)
\end{equation}
valid for $-R\log\epsilon^{-1}<\Im\tau<-R$ where $\Delta(\varphi,z)=\mathbf{\Gamma}^{u}_\epsilon(\varphi,z)-\mathbf{\Gamma}^{s}_\epsilon(\varphi,z)$. 

In order to derive an asymptotic formula for the homoclinic invariant, we consider an auxiliary function
defined by
\begin{equation*}
\Theta(\varphi,z)=\Omega\left(\Delta(\varphi,z),\partial_\varphi\mathbf{\Gamma}^{s}_\epsilon(\varphi,z)\right)\,,
\end{equation*}
where $\Omega$ is the standard symplectic form.
The homoclinic invariant of the primary homoclinic orbit is defined by \eqref{Eq:hominv}
which takes the form
\begin{equation}\label{4:homoclinicinvariant}
\omega=\Omega\bigl(\partial_\varphi \mathbf{\Gamma}^{u}_\epsilon(0,0),
\partial_\varphi \mathbf{\Gamma}^{s}_\epsilon(0,0)\bigr)\,.
\end{equation}
Differentiating the definition of $\Theta$ at the origin and taking into account that $\Delta(0,0)=0$
we get the relation:
$$
\omega=\partial_\varphi\Theta(0,0).
$$ 

Thus, we only need to estimate the function $\Theta$ and its derivative. Considering higher approximations of $u^{\pm}_\epsilon$ in \eqref{Eqs:upmapp} it is possible to improve the estimate in \eqref{Eq:est1}. In \cite{JP10} it is proved that in a neighbourhood of the point $\tau=-i\log(\epsilon^{-1})$ the following estimate holds:
\begin{equation}\label{Eq:vareq}
\Delta(\varphi,\tfrac{\beta_\epsilon}{\alpha_\epsilon}\tau+i\tfrac{\pi}{2})=-\Delta_0(\varphi,\tau)+O(\epsilon^2)
\end{equation}
which leads to
\begin{equation}\label{Eq:bound}
\begin{split}
\Theta(\varphi,z)&=-\theta_0(\varphi,\tau)+O(\epsilon^2)
=-e^{-i(\tau-\varphi)}\Theta_0(\kappa)+O(\epsilon^2)\,,
\end{split}
\end{equation}

Now note that the function $\Theta$ satisfies the following equation,
\begin{equation}\label{3:ThetaDepsilon}
(\alpha_\epsilon\partial_\varphi+\beta_\epsilon\partial_z) \Theta=\Omega(F(\Delta),\partial_\varphi\mathbf{\Gamma}^u_\epsilon),
\end{equation}
where $F(\Delta)=X_{H_{\epsilon}}(\mathbf{\Gamma}^u_\epsilon+\Delta)-X_{H_{\epsilon}}(\mathbf{\Gamma}^u_\epsilon)-DX_{H_{\epsilon}}(\mathbf{\Gamma}^u_\epsilon)\Delta$. As $F(\Delta)$ is of second order in $\Delta$ then $\Theta$ approximately satisfies the homogeneous equation $(\alpha_\epsilon\partial_\varphi+\beta_\epsilon\partial_z) u=0$ with an error of the order of $O(\left|\Delta(\varphi,z)\right|^2)$. 
Taking into account that the splitting of separatrices is rather small, we continue
our arguments neglecting this error. Then there is a function $f$ such that
$$
\Theta(\varphi,z)=f(\alpha_\epsilon  z-\beta_\epsilon\varphi)
$$
inside the domain of $\Theta$, which implies that $f$ can be extended by periodicity
onto the strip $|\Im(z)|<\tfrac\pi2-R\delta$. 
We expand the function $f$ into Fourier series, i.e.,
$$
\Theta(\varphi,z)=\sum_{k\in\mathbb{Z}}f_ke^{ik (\tfrac{\alpha_\epsilon}{\beta_\epsilon}
z-\varphi)}\,.
$$
The coefficients of the series can be expressed in terms of Fourier integrals:
\begin{eqnarray}
f_{k}&=&    \frac{\alpha_\epsilon}{2\pi\beta_\epsilon}
\int_0^{\tfrac{2\pi\beta_\epsilon}{\alpha_\epsilon}} e^{-ik\tfrac{\alpha_\epsilon}{\beta_\epsilon}z}\Theta(0,z)dz\,.
\end{eqnarray}
Following the common procedure of Fourier Analysis, we 
shift the contour of integration to $\Im z=\frac\pi2-\tfrac{\beta_\epsilon}{\alpha_\epsilon}\log\epsilon^{-1}$, 
change the variable to \eqref{4:ztaurelation} and use the estimate \eqref{Eq:bound} to get
\begin{eqnarray}
f_{-1}&=& -e^{-\tfrac{\pi\alpha_\epsilon}{2\beta_\epsilon}}\left(\Theta_0(\kappa)+O(\epsilon)\right)\,,
\end{eqnarray}
$f_1=\overline{f_{-1}}$ and there is a positive constant $C$ such that
$$
|f_k|\le C\epsilon^{2-|k|}e^{-|k|\tfrac{\pi\alpha_\epsilon}{2\beta_\epsilon}}\qquad
\mbox{for $|k|\ge2$}.
$$
Substituting these estimates into the Fourier series we get that 
for real values of $\varphi,z$
\begin{equation}\begin{split}
\Theta(\varphi,z)&=-2e^{-\frac{\pi\alpha_\epsilon}{2\beta_\epsilon}}\left|\Theta_0\right|\cos\left(\frac{\alpha_\epsilon}{\beta_\epsilon}z
-\varphi-\arg(\Theta_0)\right)+O(e^{-\frac{\pi\alpha_\epsilon}{2\beta_\epsilon}}\epsilon)\,,\\
\partial_\varphi\Theta(\varphi,z)&=-2e^{-\frac{\pi\alpha_\epsilon}{2\beta_\epsilon}}
\left|\Theta_0\right|\sin\left(\frac{\alpha_\epsilon}{\beta_\epsilon}z-\varphi-\arg(\Theta_0)\right)
+O(e^{-\frac{\pi\alpha_\epsilon}{2\beta_\epsilon}}\epsilon)\,.\\
\end{split}
\end{equation}
Since $\Theta(0,0)=0$ for all $\epsilon$ then $\arg(\Theta_0)=\pm\frac{\pi}{2}$, i.e., the Stokes constant is a purely imaginary number
and equation \eqref{4:Asymptoticformulahomoclinic} follows directly.

We note that the integrability of the normal form allows us to
repeat the arguments with more accurate approximations of the separatrices,
the result of this consideration leads to the conjecture that
\begin{equation}\label{4:AsymptoticExpansionHomoclinic}
\omega(\epsilon)\asymp e^{-\frac{\pi\alpha_\epsilon}{2\beta_\epsilon}}\sum_{k\geq 0}\omega_k \epsilon^k
\end{equation}
where $\omega_0=2\mathrm{Im}(\Theta_0(\kappa))$.

\section{Computation of the Stokes constant\label{Se:num}}

Since the arguments involved in the derivation of the asymptotic formula are not rigorous, we have developed numerical methods to check the validity of our results. 
The procedure is based on comparison of two different methods for evaluation of the
Stokes constants. The first method relies on the definition \eqref{4:SC} 
and involves the GSHE with $\epsilon=0$ only. The second method evaluates
the homoclinic invariant for $\varepsilon\ne0$ and
relies on the validity of the asymptotic expansion \eqref{4:AsymptoticExpansionHomoclinic} 
to extrapolate the values of the (normalised) homoclinic invariant towards $\varepsilon=0$
in order to get $\omega_0$.

\subsection{A method for the computation of the Stokes constant}

Let us describe the first method for computing the Stokes constant. 
We set $\tau=-i\sigma$ for $\sigma>0$, $\varphi=0$ and rewrite equation~\eqref{Eq:theta0} in the form:
\begin{equation}\label{4:theta0}
\Theta_0=\theta_0(0,-i\sigma)e^{\sigma}+O\left(e^{-(1-\epsilon_0)\sigma}\right).
\end{equation}
Then we proceed as follows. 
\begin{enumerate}
\item The first step is to construct a good approximation of stable and unstable manifolds. 
This approximation is given by a finite sum of the unique formal separatrix $\hat{u}_0$ defined in section \ref{Se:Stokes}. Given $N\geq1$ and the formal separatrix $\hat{u}_0$ we can use the relations \eqref{4:changeofvariables} to define,
\begin{equation*}
\mathbf{\Gamma}_N(\varphi,\tau):=\sum_{k=1}^N \Gamma_k(\varphi)\tau^{-k}\,,
\end{equation*} 
where
\begin{equation*}
\Gamma_k(\varphi)=\sum_{j=-k}^{k}\Gamma_{k,j}e^{ji\varphi}\text{ with } \Gamma_{k,j}\in \mathbb{C}^4,
\end{equation*}
such that $\mathbf{\Gamma}_N$ approximates the parameterizations $\mathbf{\Gamma}^{\pm}_0$ in the following sense
\begin{equation*}
\mathbf{\Gamma}^{\pm}_0(\varphi,z)-\mathbf{\Gamma}_N(\varphi,\tau)=O(\tau^{-N-1})\,.
\end{equation*}
The natural number $N$ can be chosen using the \textsl{astronomers recipe}. 
It simply chooses $N$ such that for fixed $\tau$ and $\varphi$ it minimizes $\left|\Gamma_{N+1}(\varphi)\tau^{-N-1}\right|$, 
that is, the least term of the formal series $\sum_{k\geq1}\Gamma_k(\varphi)\tau^{-k}$ (see Figure \ref{4:figAst2}).
\begin{figure}
\centering
      \includegraphics[width=2in]{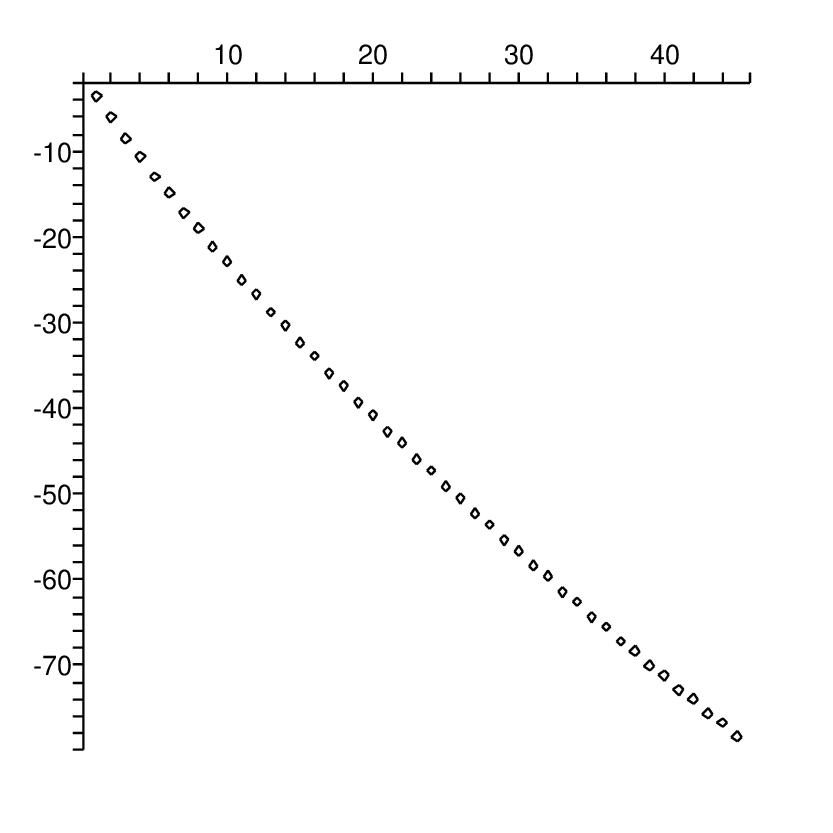}%
\caption{Graph of $\log_{10}\left(\frac{\max_{j}\left|\Gamma_{k,j}\right|}{(350\pi)^k}\right)$}.
\label{4:figAst2}
%\end{minipage}
\end{figure}
\item A point on the unstable manifold (resp. stable manifold) can be represented in the coordinates $(\varphi,\tau)$. 
In order to obtain a point close to the unstable manifold we fix a positive real number $\sigma \in \mathbb{R}^{+}$ 
and a sufficiently large $d \in \mathbb{R}^{+}$ and define $z^{-}_0=\mathbf{\Gamma}_N(-d,-i\sigma-d)$ and 
a tangent vector $v^{-}_0=\partial_\varphi\mathbf{\Gamma}_N(-d,-i\sigma-d)$. Analogously, for the stable manifold 
we define $z^{+}_0=\mathbf{\Gamma}_N(d,-i\sigma+d)$ and $v^{+}_0=\partial_\varphi\mathbf{\Gamma}_N(d,-i\sigma+d)$. 

\item The next step is to measure the difference of stable and unstable manifold at the point $(\varphi,\tau)=(0,-i\sigma)$. 
Taking into account the periodicity in $\varphi$ we set $d$ equal to a multiple to $2\pi$ and integrate numerically the system,
\begin{equation}\label{4:S1}\begin{array}{l}
z'=X_{H_0}(z)\\
v'=DX_{H_0}(z)v\\
\end{array}
\end{equation}
forward in time with $t \in [0,d]$ and initial conditions $z^-(0)=z_0^{-},v^-(0)=v_0^{-}$ 
and then backward in time with $t \in [-d,0]$ and initial conditions $z^+(0)=z_0^{+},v^+(0)=v_0^{+}$. 
\item Finally we evaluate,
\begin{equation}\label{4:theta}
\hat{\Theta}(\sigma)=\Omega(z^+(-d)-z^-(d),v^-(d))e^{\sigma}
\end{equation}
\end{enumerate}
\begin{remark}
The stable and unstable manifolds have the same asymptotic expansion, hence the difference $z^+(-d)-z^-(d)$ is exponentially small,
i.e. comparable with $e^{\sigma}$. Thus the system \eqref{4:S1} has to be integrated with great accuracy. 
In the case of GSHE an excellent integrator can be constructed using a high order Taylor series method. 
\end{remark}

\subsection{Numerical results}
In all current computations we have used a Taylor series method, which is incorporated in the Maple Software, 
to integrate the equations of motion \eqref{4:S1}. The method uses an adaptive step procedure controlled 
by a local error tolerance which was set to $10^{-D}$, where $D$ is the number of significant digits used 
in the computations. The order of the method has been automatically defined using the formula $\max(22,\left\lfloor 1.5D\right\rfloor)$.

Having fixed $\kappa=2$ (which we recall to be one of the parameters of the original equation \eqref{4:GSHE}) we have computed the first 45 coefficients of the formal separatrix $\hat{u}_0$ 
with 60 digits precision. The error committed by the approximation $\mathbf{\Gamma}_N$ is approximately of the order of the first missing term.
\begin{figure}[tbp]
  \begin{center}
    \includegraphics[width=4in]{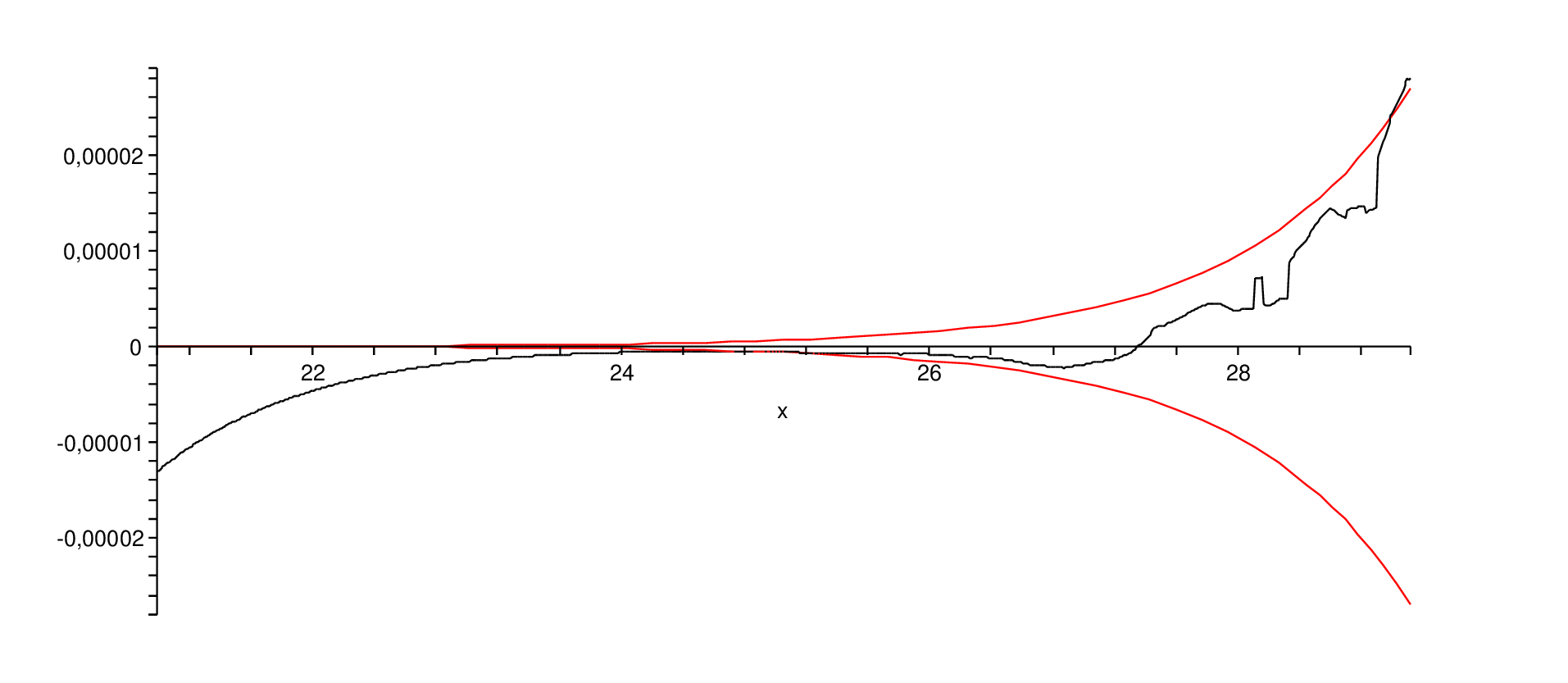}
  \end{center}
\begin{center}
    \includegraphics[width=4in]{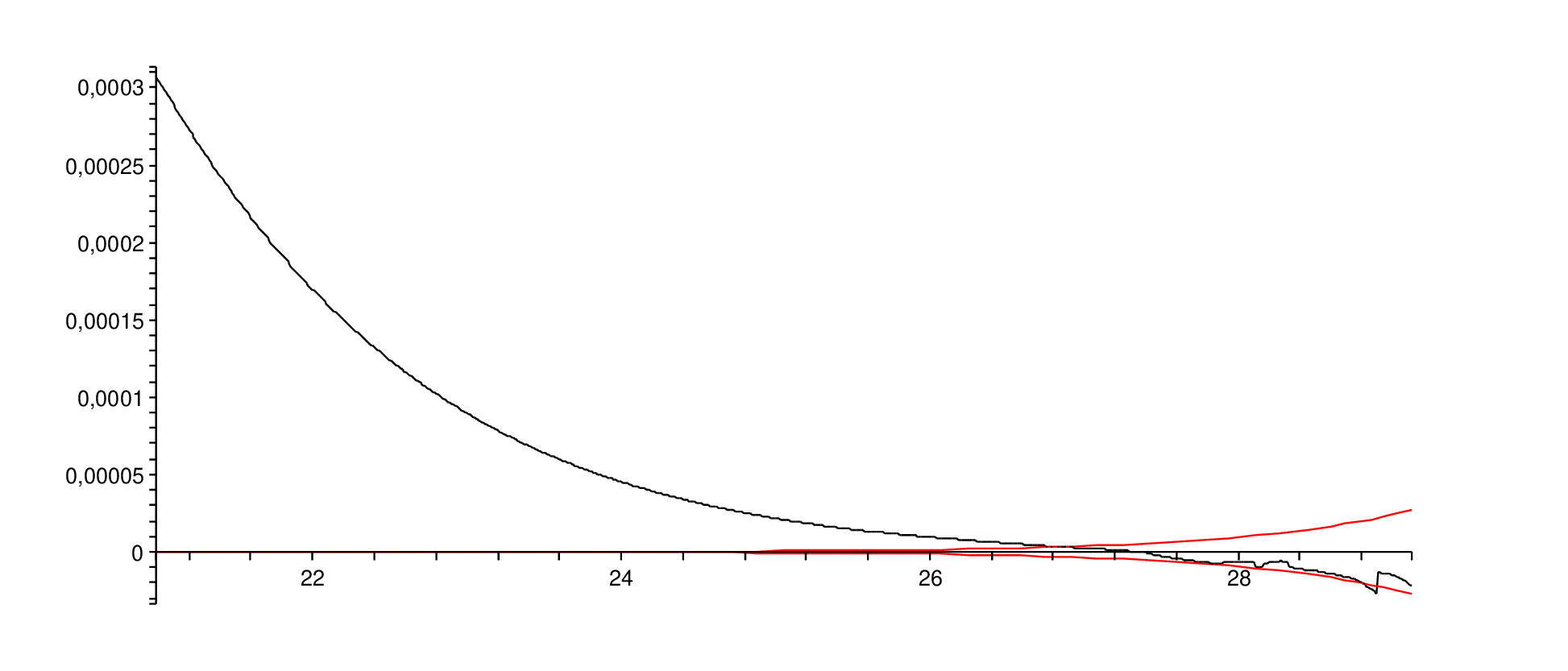}
  \end{center}
  \caption{\small The top figure represents the graph of the function 
$\mathrm{Im}(\hat{\Theta}(\sigma))e^{\sigma}-10.472161956944$ and the bottom figure represents 
the graph of the function $\mathrm{Re}(\hat{\Theta}(\sigma))e^{\sigma}$. When $\sigma$ 
is around 25 the rounding errors become visible and the convergence stops. 
The dashed curves represent the magnitude of the rounding errors.}
  \label{4:fig1}
\end{figure}

Using double precision (16 digits) we have integrated numerically the equations \eqref{4:S1} 
to obtain $\hat{\Theta}(\sigma)$ for values of $\sigma$ uniformly distributed in the interval 
$\left[20,28.89\right]$. The initial conditions were computed using $d=350\pi$ and the first
 9 terms of $\mathbf{\Gamma}_N$. The results are depicted in Figure \ref{4:fig1}. 
The expected errors are bounded by the dashed curves. This implies in particular that the method is numerically 
stable, that is, the propagation errors due to integration do not increase drastically.
There are several sources of errors that affect the accuracy of the computation of the Stokes constant, namely:
\begin{itemize}
	\item Approximation of stable and unstable manifolds given by the function $\mathbf{\Gamma}_N$;
	\item Errors due to the numerical integration;
	\item Rounding errors.
\end{itemize}
The first and the second source of errors can be made small compared to the rounding errors, which can be roughly estimated by,
\begin{equation}\label{4:roundoff}
\frac{C}{\sigma^2}10^{-D}e^{\sigma},
\end{equation}
where $D$ is the number of digits used in the computations and $C$ is a real positive constant 
which reflects the propagation of rounding errors. Using this estimate we have provided bounds 
for the rounding errors which can be observed in Figure \ref{4:fig1}. The constant $C$ can 
be estimated by fitting the function \eqref{4:roundoff} to the points $\left|\hat{\Theta}(\sigma)\right|$ for $\sigma\geq 25$. 
Using the method of least squares we have concluded that $C$ is approximately $38.5$.

With double arithmetic precision the method previously described allows the computation of 7 to 8 correct digits 
of the Stokes constant $\Theta_0$. In fact the rounding errors in computing $\hat{\Theta}(\sigma)$ from 
formula \eqref{4:theta} grow accordingly to \eqref{4:roundoff} whereas the neglected terms of the formula 
\eqref{4:theta0} decrease like $C_1e^{-\sigma}$, where $C_1$ is some positive constant. Hence the optimal 
is attained when both contributions are of the same order. The constant $C_1$ can be estimated by fitting 
the function $C_0+C_1e^{-\sigma}$ to the points $\left|\hat{\Theta}(\sigma)\right|$ for $\sigma\leq 24$.
 Using the method of least squares we have obtained that $C_1$ is approximately $17305.75$.
Using this information we can determine the value $\sigma^*$ where both contributions are essentially of the same order. 
This means that $\sigma^*$ must satisfy the equation,
\begin{equation*}
(e^{-\sigma})^2=\frac{C}{{\sigma}^2\,C_1}10^{-D}
\end{equation*}
which implies,
\begin{equation*}
\left|\Theta_0-\hat{\Theta}(\sigma^*)\right|\approx \frac{816}{\sigma^*}10^{-\frac{D}{2}}
\end{equation*}
\begin{table*}[tp]
	\tiny
	\begin{tabular}{|c|c|c|l|}
	\hline
	$D$&$\sigma^*$&$\mathrm{Re}(\hat{\Theta}(\sigma^*))$&$\mathrm{Im}(\hat{\Theta}(\sigma^*))$\\
	\hline
	16&24.68&2.7e-05&\textbf{10.472161}43901571\\
	20&29.46&7.8e-07&\textbf{10.47216195}3423286113\\
	24&34.21&1.6e-08&\textbf{10.4721619569}069446924024\\
	28&38.95&3.1e-10&\textbf{10.472161956944}13924682820786\\
	32&43.67&5.3e-12&\textbf{10.47216195694439}6725504278408504\\
	36&48.37&8.5e-14&\textbf{10.4721619569443983}419527788851129556\\
	40&53.07&1.2e-15&\textbf{10.472161956944398358}12989263311456886391\\
	44&57.76&1.8e-17&\textbf{10.47216195694439835828}4180684468467819622191\\
	48&62.45&2.6e-19&\textbf{10.4721619569443983582855}084356725900717201861670\\
	52&67.12&3.5e-21&\textbf{10.472161956944398358285521}30242825730920048239485015\\
	56&71.80&4.7e-23&\textbf{10.47216195694439835828552143}0879142372532568396894067732\\
	60&76.46&6.2e-25&\textbf{10.4721619569443983582855214320}209319731283197852962601326570\\
	64&81.13&8.0e-27&\textbf{10.472161956944398358285521432031}66495538939445255794702026972749\\
	68&85.79&1.0e-28&\textbf{10.47216195694439835828552143203190}0047829633854060398152634432422925\\
	\hline
	\end{tabular}
	\vspace{10pt}
		\caption{Stokes constant evaluated at the optimum $\sigma^*$ for different computer precisions. 
In the computations we have used $d=350\pi$ and $N=40$}
	\label{4:tabStokesConstant}
\end{table*}

In this way it is possible to obtain 8 correct digits for the Stokes constant using only double precision.
In Table \ref{4:tabStokesConstant} we have listed the values of $\hat{\Theta}(\sigma^{*})$ evaluated 
at the optimum $\sigma^*$ for higher computer precisions. The digits in bold correspond 
to correct digits of the Stokes constant. We also note that the numerics suggest that $\Theta_0$ 
is pure imaginary which agrees with our prediction.

Finally, let us mention that in the process of computing the Stokes constant we have made several 
choices for the parameters. Namely, the number of terms $N$ used to compute $\mathbf{\Gamma}_N$ 
and the parameter $d$ which were used in computing the initial conditions of step (ii) of the numerical scheme. 
In fact the results are independent of these particular choices and Table \ref{4:tabrobust} demonstrates 
the robustness of the numerical method.
\begin{table*}[tp]
	\begin{tabular}{|c|c|c|c|}
	\hline
	$d \backslash N$&10&20&30\\
	\hline
	100$\pi$&10.47216215179386&10.47216215183208&10.47216215181955\\
	\hline
	150$\pi$&10.47216131335742&10.47216131335746&10.47216131335772\\
	\hline
	200$\pi$&10.47216144775669&10.47216144775671&10.47216144775682\\
	\hline
	250$\pi$&10.47216149546998&10.47216149546998&10.47216149547027\\
	\hline
	300$\pi$&10.47216132022817&10.47216132022820&10.47216132022773\\
	\hline
	350$\pi$&10.47216138600882&10.47216138600883&10.47216138600868\\
	\hline
	\end{tabular}
	\vspace{10pt}
		\caption{Comparison of the value of $\mathrm{Im}(\hat{\Theta}(25))$ for different values of parameters $N$ and $d$.}
		\label{4:tabrobust}
\end{table*}

\section{High precision computations of an asymptotic expansion for the homoclinic invariant\label{Se:num3}}
In this section we present a numerical method for the computation of the homoclinic invariant as defined 
in \eqref{4:homoclinicinvariant} for the Swift-Hohenberg equation with $\kappa=2$ and $\epsilon<0$. 
Moreover we investigate from a numerical point of view the validity of the asymptotic expansion 
\eqref{4:AsymptoticExpansionHomoclinic} for the homoclinic invariant. This section follows the ideas of \cite{GS:08}
originally developed for the study of exponentially small phenomena for area-preserving maps.

%Recall the definition of the homoclinic invariant,
%\begin{equation}\label{4:homoclinicinvariant2}
%\omega(\epsilon)=\left.\Omega(\partial_\phi \mathbf{\Gamma}^{u}_\epsilon(\alpha_\epsilon \phi,\beta_\epsilon s),
%\partial_\phi \mathbf{\Gamma}^{s}_\epsilon(\alpha_\epsilon \phi,\beta_\epsilon s))\right|_{(\phi,s)=(0,0)}
%\end{equation}
%where $\mathbf{\Gamma}^{s,u}_\epsilon$ are solutions of the following linear PDE,
%\begin{equation}\label{4:PDE_formal}
%(\alpha_\epsilon\partial_\varphi + \beta_\epsilon\partial_z)\mathbf{\Gamma}= X_{H_\epsilon}(\mathbf{\Gamma})
%\end{equation}
In order to compute the homoclinic invariant \eqref{Eq:hominv} 
we need to compute two tangent vectors at the symmetric homoclinic point $\mathbf{\Gamma}^{s}_\epsilon(0,0)$. 
Using the fact that the system is reversible we can obtain the stable tangent vector 
$\partial_\varphi \mathbf{\Gamma}^{s}_\epsilon$ by applying the reversor to the unstable 
tangent vector $\partial_\varphi \mathbf{\Gamma}^{u}_\epsilon$. The unstable tangent vector 
$\partial_\varphi \mathbf{\Gamma}^{u}_\epsilon$ lives in the tangent plane of the unstable 
manifold at the symmetric homoclinic orbit. Thus an easy way to compute this tangent vector 
is to approximate the primary homoclinic orbit near the equilibrium point by the following expansion,
\begin{equation}\label{4:FormalExpansion}
	\mathbf{\Gamma}_{\epsilon,N}^{u}(\varphi,z)=\sum_{k=1}^{N}e^{k z}\left(\mathbf{c}_k(\epsilon)
+\sum_{j\geq1}^k \mathbf{a}_{k,j}(\epsilon)\cos(j\varphi) + \mathbf{b}_{k,j}(\epsilon)\sin(j\varphi)\right)
\end{equation}
and then use the variational equations,
\begin{equation}\label{4:equationsofmotions1}
\begin{split}
\mathbf{x}'&=X_{H_\epsilon}(\mathbf{x})\\
\mathbf{v}'&=DX_{H_\epsilon}(\mathbf{x})\mathbf{v}\\
\end{split}
\end{equation}
to transport the tangent vector $\partial_\varphi\mathbf{\Gamma}_{\epsilon,N}^{u}$ along 
the primary homoclinic orbit until it hits the symmetric plane $\mathrm{Fix}(S)$ defined by $\left\{q_2=0,p_1=0\right\}$. 
Let us present the details of the method.
\subsection{A method for the computation of the homoclinic invariant}
\begin{enumerate}
	\item The first step is to determine the coefficients of \eqref{4:FormalExpansion}. To that end we take a new expansion,
\begin{equation*}
u_N(\varphi,z)=\sum_{k=1}^{N}e^{k z}\left(c_k(\epsilon)+\sum_{j\geq1}^k a_{k,j}(\epsilon)\cos(j\varphi) + b_{k,j}(\epsilon)\sin(j\varphi)\right)
\end{equation*}
and substitute into the equation,
\begin{equation}\label{4:eqSHnormal}
((\alpha_\epsilon\partial_\varphi+\beta_\epsilon\partial_z)^2+1)^2\,u=\epsilon u+2u^2-u^3
\end{equation} 
and collect the terms of the same order in $e^{kz}$. In this way it is possible to determine coefficients 
$c_k$, $a_{k,j}$ and $b_{k,j}$. It is not difficult to see that the coefficients $a_{1,1}$ and $b_{1,1}$ 
satisfy no relations and that all other coefficients depend from these two. So we define,
\begin{equation*}
a_{1,1}=r_0\cos(\psi_0)\ \ \mathrm{and} \ \ b_{1,1}=r_0\sin(\psi_0)
\end{equation*}
Now recall that the first component of $\mathbf{\Gamma}^{u}_\epsilon$ solves equation \eqref{4:eqSHnormal} 
and due to the asymptotic behavior \eqref{4:leadingorderGamma_epsilon} we conclude that for $z<<0$ and $\delta<<1$ 
it is approximately,
\begin{equation}\label{4:leadingordersGamma}
e^{z}\left(-\frac{2\delta}{\sqrt{\eta}}\cos(\varphi)+\frac{\delta^2}{\sqrt{\eta}}\left(1+\frac{2\mu}{\eta}\right)\sin(\varphi)\right)+O(e^{2z})
\end{equation}
where $\epsilon=-4\delta^2$. Next we \textsl{"match"} the leading order of $u_N(\phi,s)$ with the expression 
\eqref{4:leadingordersGamma} and conclude that $\psi_0$ and $r_0$ must satisfy,
\begin{equation}\label{4:choiceparameters}
\begin{split}
\psi_0&=\arctan\left(-\left(1+\frac{2\mu}{\eta}\right)\frac{\delta}{2}\right)\\
r_0&=\frac{2\delta}{\sqrt{\eta}}\sqrt{1+\left(1+\frac{2\mu}{\eta}\right)^2\frac{\delta^2}{4}}
\end{split}
\end{equation}
Taking into account \eqref{4:changeofvariables} we reconstruct $\mathbf{\Gamma}_{\epsilon,N}^{u}$ from $u_N$ 
and due to the \textsl{"matching"} \eqref{4:choiceparameters} we have,
\begin{equation*}
\mathbf{\Gamma}_{\epsilon}^{u}(t,t)\approx \mathbf{\Gamma}_{\epsilon,N}^{u}(t,t), \ \mathrm{as} \ \ t\rightarrow -\infty,\ \delta\rightarrow0.
\end{equation*}
That is, for small values of $\delta$, the expansion $\mathbf{\Gamma}_{\epsilon,N}^{u}$ provides 
a good approximation of the primary homoclinic orbit near the equilibrium point.
\item The second step is to improve the accuracy of the approximation of the symmetric homoclinic point, 
provided by $\mathbf{\Gamma}_{\epsilon,N}^{u}$. Given small $\delta$ and sufficiently large $T_0>0$ we want to determine $(T,\psi)$ such that,
\begin{align*}
\mathbf{x}'&=X_{H_\epsilon}(\mathbf{x}),&\mathbf{x}(0;\psi)&=\mathbf{\Gamma}_{\epsilon,N}^{u}(-\alpha_\epsilon T_0,-\beta_\epsilon T_0;\psi)
\end{align*}
subject to,
\begin{equation}\label{4:reversibilitycondition}
\mathbf{x}(T;\psi) \in \mathrm{Fix}(S)
\end{equation}
This problem can be solved using Newton method. Starting from $(T_0,\psi_0)$ we obtain a sequence of points $(T_i,\psi_i)$,
\begin{equation}\label{4:newtonmethod}
\begin{pmatrix}
T_{i+1}\\
\psi_{i+1}
\end{pmatrix}=\begin{pmatrix}
T_{i}\\
\psi_{i}
\end{pmatrix}-\begin{pmatrix}
\frac{\partial q_2}{\partial T}(T_i;\psi_i)&\frac{\partial q_2}{\partial \psi}(T_i;\psi_i)\\
\frac{\partial p_1}{\partial T}(T_i;\psi_i)&\frac{\partial p_1}{\partial \psi}(T_i;\psi_i)
\end{pmatrix}^{-1}\begin{pmatrix}
q_2(T_i;\psi_i)\\
p_1(T_i;\psi_i)
\end{pmatrix}
\end{equation}
that converges to a limit $(T_{*},\psi_{*})$ such that $\mathbf{x}(T_*;\psi_*) \in \mathrm{Fix}(S)$, 
provided $(T_0,\psi_0)$ is sufficiently close to $(T_*,\psi_*)$ (see \cite{ARCAS:93}). The derivatives 
in \eqref{4:newtonmethod} can be computed using the variational equations along the orbit $\mathbf{x}(t;\psi)$. 
Later we will see that the formulae \eqref{4:choiceparameters} provide sufficiently accurate initial guesses 
yielding the convergence of the Newton method. 

\item Having obtained in the previous step an accurate approximation of the symmetric homoclinic point, 
the last step is to integrate numerically the system,
\begin{align*}
\mathbf{x}'&=X_{H_\epsilon}(\mathbf{x}),&\mathbf{x}(0;\psi)&=
\mathbf{\Gamma}_{\epsilon,N}^{u}(-\alpha_\epsilon T_0,-\beta_\epsilon T_0;\psi_{*})\\
\mathbf{v}'&=DX_{H_\epsilon}(\mathbf{x})\mathbf{v},& \mathbf{v}(0;\psi)
&=\alpha_\epsilon\partial_\varphi \mathbf{\Gamma}_{\epsilon,N}^{u}(-\alpha_\epsilon T_0,-\beta_\epsilon T_0;\psi_{*})\\
\end{align*}
and evaluate the homoclinic invariant,
\begin{equation*}
\hat{\omega}=\Omega(\mathbf{v}(T_{*},\psi_{*}),S(\mathbf{v}(T_{*},\psi_{*})))
\end{equation*}
\end{enumerate}
\subsection{Numerical results}
\begin{figure}
  \begin{center}
    \includegraphics[width=4in]{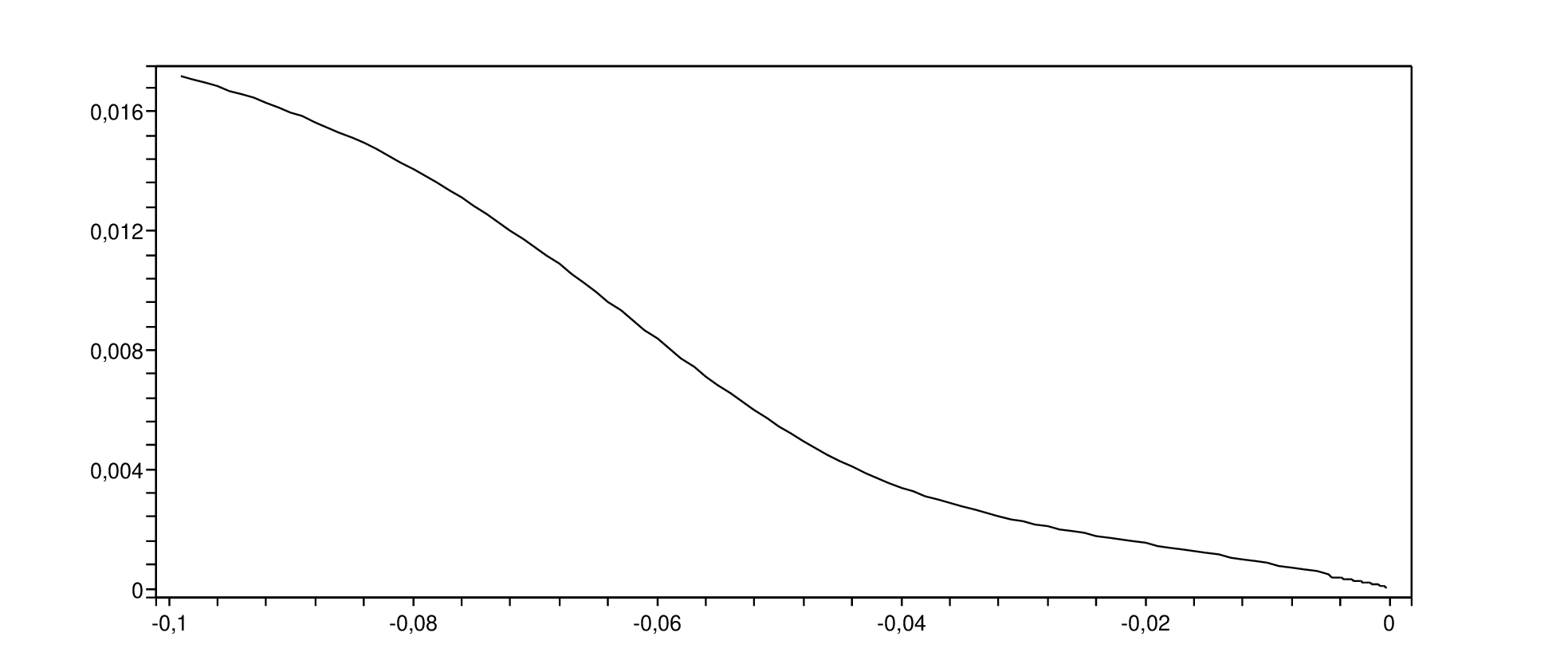}
  \end{center}
  \caption{\small Relative error of $(T_0,\psi_0)$ depending on $\epsilon\in\mathcal{I}$}.
  \label{4:figrelativeerrornewtonmethod}
\end{figure}

We have considered a finite set $\mathcal{I}$ consisting of points in the interval 
$\epsilon \in [-\frac{1}{10},-\frac{1}{1000}]$ and computed the homoclinic invariant 
for those points using the method previously described. For all points in $\mathcal{I}$ 
the magnitude of homoclinic invariant ranges from $10^{-5}$ to $10^{-45}$. 
Thus, in all numerical integrations we have used a high order Taylor method which 
allows to perform the numerical integration with very high precision. We have computed 
the coefficients of the expansion \eqref{4:FormalExpansion} up to $N=5$ and for each 
$\epsilon\in \mathcal{I}$ we have chosen $T_0$ sufficiently large so that 
$\mathbf{\Gamma}_{\epsilon,N}^{u}(-\alpha_\epsilon T_0,-\beta_\epsilon T_0)$ approximates 
the unstable manifold within the required precision. The initial point $(T_0,\psi_0)$ 
used in Newton method proved to be very close to $(T_*,\psi_*)$ and its relative 
error can be observed in Figure \ref{4:figrelativeerrornewtonmethod}.
After computing the homoclinic invariant we have normalized it using the formula,
\begin{equation*}
\bar{\omega}(\epsilon)=\frac{\omega(\epsilon)}{2}e^{\frac{\pi\alpha_\epsilon}{2\beta_\epsilon}}
\end{equation*}
\begin{figure}
  \begin{center}
    \includegraphics[width=3in]{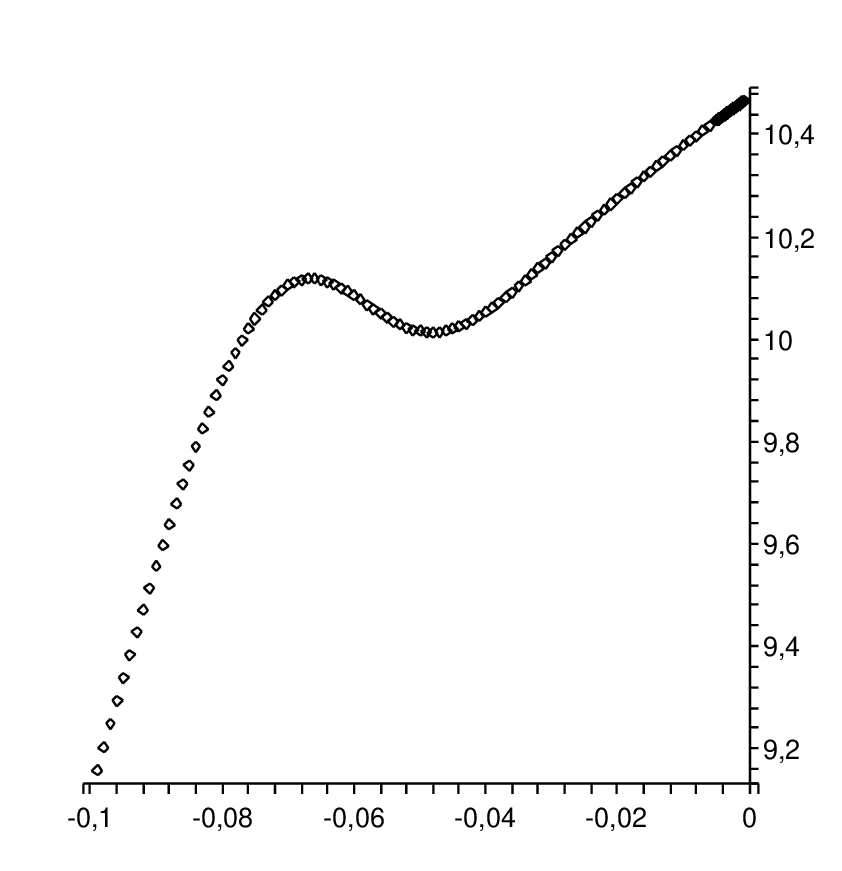}
  \end{center}
  \caption{\small Graph of the function $\bar{\omega}(\epsilon)$}.
  \label{4:fighomoclinic}
\end{figure}
The behaviour of the function $\bar{\omega}(\epsilon)$ can be observed in Figure \ref{4:fighomoclinic}. 
It possible to see that it is approaching the value of the Stokes constant computed in the previous section. 
Moreover, it is aproaching this value in a linear fashion, supporting the validity 
of the asymptotic formula \eqref{4:Asymptoticformulahomoclinic}. 
Taking into account the asymptotic expansion for $\omega(\epsilon)$ we investigate 
the validity of the following asymptotic expansion for $\bar{\omega}(\epsilon)$,
\begin{equation}\label{4:asymptoticexpansion2}
\bar{\omega}(\epsilon)\asymp\sum_{k\geq0}\bar{\omega}_k\epsilon^k
\end{equation}%

\begin{table*}[tp]
\center
\scriptsize
	\begin{tabular}{|c|l|l|l|}
	\hline
	&$\bar{\omega}_0$&$\bar{\omega}_1$&$\bar{\omega}_2$\\
	\hline
 5&10.47216195694& 8.979943127&- 42.60110\\6&10.472161956944&
 8.979943127&- 42.601100\\7&
 10.4721619569443& 8.9799431275&- 42.6011004\\8& 10.47216195694439& 8.97994312752&- 42.60110043\\9& 10.472161956944398&
 8.9799431275209&- 42.601100432
\\10& 10.4721619569443983& 8.9799431275210&-
 42.601100432\\11&
 10.4721619569443983& 8.9799431275210&- 42.601100432\\12& 10.4721619569443983&
 8.9799431275210&- 42.6011004327\\
	\hline
	&$\bar{\omega}_3$&$\bar{\omega}_4$&$\bar{\omega}_5$\\
	\hline
	5& 152.88&- 774.4& 3.8$\times 10^3$\\
	6& 152.888&- 774.2& 3.8$\times 10^3$\\
  7& 152.887&- 774.40& 3.80$\times 10^3$\\
  8& 152.88795&- 774.39& 3.814$\times 10^3$\\
  9& 152.88795&- 774.394& 3.813$\times 10^3$\\
  10&152.887958&- 774.3944& 3.8138$\times 10^3$\\
  11& 152.887958&- 774.3944& 3.813$\times 10^3$\\
  12& 152.887958&- 774.3944& 3.813$\times 10^3$ \\
  \hline
	\end{tabular}
	\vspace{10pt}
		\caption{Coefficients of the estimated polynomials for different subsets of $\mathcal{P}$ and different degrees.}
		\label{4:Tableestimatedcoefficients}
\end{table*}

To that end, we have taken 14 points evenly spaced in the interval $[-2.7\times10^{-3},-1.4\times10^{-3}]$ 
and computed the corresponding normalized homoclinic invariant with more than 40 correct digits. 
Let us denote this set of homoclinic invariants by $\mathcal{P}$.
Then, in order to get the first few coefficients of the asymptotic expansion \eqref{4:asymptoticexpansion2} 
we have fitted a partial sum of the asymptotic expansion to the points of $\mathcal{P}$. 
Here we have used as many points as the number of unknown coefficients.
Moreover, following \cite{GS:08} we have performed the following tests to evaluate the validity of the asymptotic expansion:
\begin{enumerate}
	\item Interpolating different partial sums to different subsets of $\mathcal{P}$ 
should give essentially the same results for the coefficients.
	\item The constant term of the interpolating polynomial should coincide with 
the value of the Stokes constant computed in the previous section.
	\item The interpolating polynomial should reasonably approximate 
$\bar{\omega}(\epsilon)$ outside the interval $[-2.7\times10^{-3},-1.4\times10^{-3}]$, 
in the sense that it agrees with the main property of an aymptotic expansion:
	\begin{equation*}
	\left|\bar{\omega}(\epsilon)-\sum_{k\geq0}^{n-1}\bar{\omega}_k\epsilon^k\right|\leq C \epsilon^n,\ \forall \epsilon \in \left[\epsilon_0,0\right)
	\end{equation*}
	for some $C>0$ and $\epsilon_0<0$.
\end{enumerate}

\begin{figure}
  \begin{center}
    \includegraphics[width=3in]{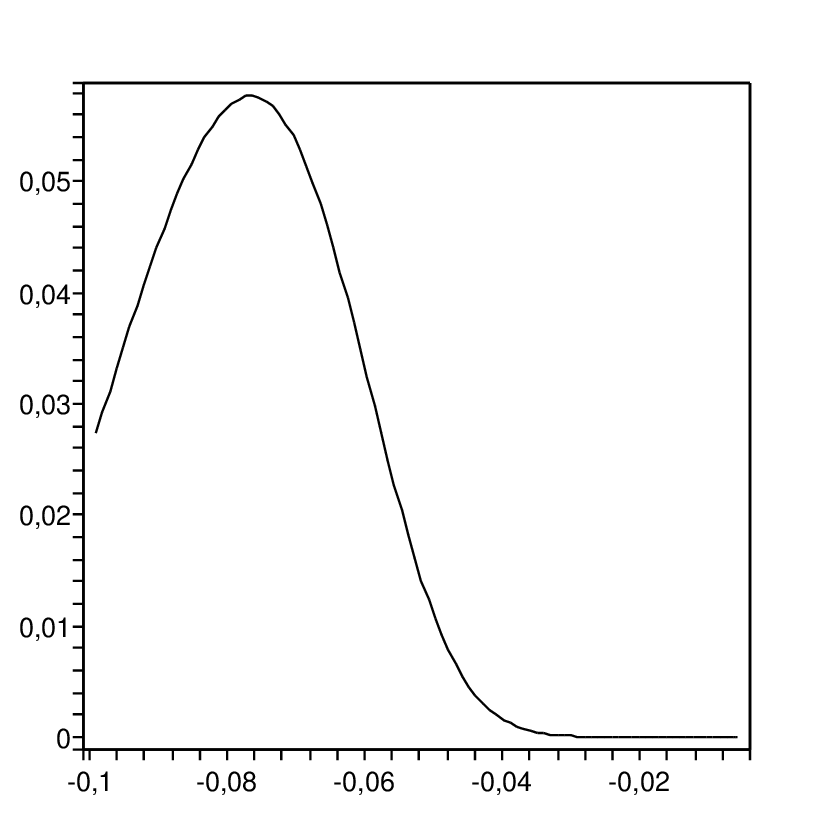}
  \end{center}
  \caption{Relative error of the asymptotic expansion of $\bar{\omega}(\epsilon)$.}
  \label{4:fig_relativeerror}
\end{figure}

For the first test we have considered all possible subsets of $\mathcal{P}$ 
having only $6$ consecutive elements and interpolated these data by polynomials 
of degree 5. Then for each coefficient, we extracted the part of the number which 
is equal to all polynomials. We have repeated this process for polynomials of degree 6 
up to degree 12. The results are summarized in Table \ref{4:Tableestimatedcoefficients}, 
where it is possible to see that there is a good agreement between the coefficients 
of the different interpolating polynomials of different subsets of $\mathcal{P}$. 
We can also infer from Table \ref{4:Tableestimatedcoefficients} that the results 
are numerically stable. Thus, we have the following estimates for the first 6 coefficients of \eqref{4:asymptoticexpansion2}:
\begin{align*}
\bar{\omega}_0&=10.4721619569443983\ldots & \bar{\omega}_1&=8.9799431275210\ldots& \bar{\omega}_2&=-42.601100432\ldots\\
\bar{\omega}_3&=152.887958\ldots          & \bar{\omega}_4&=-774.3944\ldots      & \bar{\omega}_5&=3.813\ldots\times 10^3
\end{align*}

Furthermore, it is clear that the coefficient $\bar{\omega}_0$ coincides (up to 18 digits) with the value of the Stokes constant which we recall,
\begin{equation*}
\left|\Theta_0\right|=10.47216195694439835828552143203190\ldots
\end{equation*} 
Moreover, in Figure \ref{4:fig_relativeerror} we see that the relative error of the asymptotic 
expansion does not exceed $0.06$ in the hole interval $\left[-\frac{1}{10},0\right]$. 
Thus, our numerical results provide a satisfactory numerical evidence that supports the correctness 
of the asymptotic expansion \eqref{4:AsymptoticExpansionHomoclinic}.

\begin{appendix}

\section{Transformation of GSHE to the normal form\label{Ap:A}}

In order to normalize $H_\epsilon$ up to order $5$, we have used the method 
of Lie series to determine Hamiltonians $F_i$, $i=0,\ldots,4$ which generate the following near identity canonical map,
$$
\Psi_5=\Phi_{F_0}^1\circ\Phi_{F_1}^1\circ\Phi_{F_2}^1\circ\Phi_{F_3}^1\circ\Phi_{F_4}^1\,,
$$
where
\begin{equation}\begin{split}
F_0&=\epsilon\, \left( -{\frac {5}{32}}\,{ q_1}\,{ p_1}+{\frac {3}{32}}
\,{ q_2}\,{ p_2}+\frac{1}{8}\,{ p_1}\,{ p_2} \right)\\
F_1&={\frac {7}{216}}\,\kappa\,\sqrt {2}{{ q_1}}^{2}{ p_2}+{\frac {95}{
216}}\,\kappa\,\sqrt {2}{ q_1}\,{ q_2}\,{ p_1}+{\frac {17}{72}}
\,\kappa\,\sqrt {2}{ q_1}\,{{ p_1}}^{2}+{\frac {5}{36}}\,\kappa\,
\sqrt {2}{ q_1}\,{{ p_2}}^{2}+\\
&\quad {\frac {175}{432}}\,\kappa\,\sqrt {2
}{{ q_2}}^{2}{ p_2}+\frac{1}{36}\,\kappa\,\sqrt {2}{ q_2}\,{ p_1}\,{
 p_2}-\frac{1}{12}\,\kappa\,\sqrt {2}{{ p_1}}^{2}{ p_2}-\frac{1}{18}\,\kappa\,
\sqrt {2}{{ p_2}}^{3}\\
F_2&= \left( -{\frac {517}{20736}}\,{\kappa}^{2}+{\frac {29}{512}} \right) 
{ q_1}\,{{ p_1}}^{3}+ \left( -{\frac {217}{20736}}\,{\kappa}^{2}+{
\frac {17}{512}} \right) { q_1}\,{ p_1}\,{{ p_2}}^{2}+\\
&\quad \left( {
\frac {2327}{20736}}\,{\kappa}^{2}-{\frac {31}{512}} \right) { q_2}
\,{{ p_1}}^{2}{ p_2}+ \left( -{\frac {19}{512}}+{\frac {2027}{
20736}}\,{\kappa}^{2} \right) { q_2}\,{{ p_2}}^{3}+ \\
&\quad\left( -{
\frac {5}{128}}+{\frac {7}{192}}\,{\kappa}^{2} \right) {{ p_1}}^{3}{
 p_2}+ \left( {\frac {19}{576}}\,{\kappa}^{2}-{\frac {3}{128}}
 \right) { p_1}\,{{ p_2}}^{3}\\
F_3&=\epsilon\, \left( -{\frac {143}{1152}}\,\kappa\,\sqrt {2}{{ p_1}}^{2
}{ p_2}-{\frac {167}{1728}}\,\kappa\,\sqrt {2}{{ p_2}}^{3}
 \right)\\
F_4&=-{\frac {2}{1215}}\,\sqrt {2}\kappa\, \left( 37\,{\kappa}^{2}-27
 \right) {{ p_2}}^{5}-{\frac {1}{648}}\,\sqrt {2}\kappa\, \left( -45
+52\,{\kappa}^{2} \right) {{ p_1}}^{4}{ p_2}-\\
&\quad{\frac {1}{243}}\,
\sqrt {2}\kappa\, \left( -27+34\,{\kappa}^{2} \right) {{ p_1}}^{2}{{
 p_2}}^{3}\\
\end{split}
\end{equation}
Using an algebraic manipulator it is not difficult to see that $\Psi_5$ transforms $H_\epsilon$ into the desired form.
\end{appendix}

\bibliographystyle{plain}

\section*{Acknowledgements}
The authors thank Lev Lerman for suggesting a set up of the problem
and usefull discussions.

\end{document}